\documentclass[10pt]{article}
\input epsf.tex
\usepackage{latexsym}
\usepackage{amsmath,amsthm,amsfonts,amscd,mathrsfs,pifont}  \usepackage{eucal} % Euler fonts
\usepackage{epsfig}
\usepackage{mathrsfs}
\usepackage{verbatim}     % I need the verbatim package here.

\newtheorem{theorem}{Theorem}[section]
\newtheorem{lemma}[theorem]{Lemma}
\newtheorem{proposition}[theorem]{Proposition}

\newtheorem{example}[theorem]{Example}
\newtheorem{remark}[theorem]{Remark}

\newtheorem{fact}[theorem]{Fact}

\newtheorem*{notation}{Notation}

\title{Derived Arithmetic Fuchsian Groups of Genus Two}

\date{}

\begin{document}
\title{\bf Derived Arithmetic Fuchsian Groups of Genus Two}
\author{\sc Melissa L.~Macasieb}
\maketitle

\begin{abstract} We classify all cocompact torsion-free derived arithmetic Fuchsian groups of genus two by commensurability class.  In particular, we show that there exist no such groups arising from quaternion algebras over number fields of degree greater than 5.  We also prove some results on the existence and form of maximal orders for a class of quaternion algebras related to these groups.  Using these results in conjunction with a computer program, one can determine an explicit set of generators for each derived arithmetic Fuchsian group containing a torsion-free subgroup of genus two.  We show this for a number of examples.

\end{abstract}

%\subjclass[2000]{57M50, 11R52}
%\keywords{hyperbolic orbifolds, arithmetic Fuchsian groups, maximal orders}

\newcommand\real{\mathbb R}
\newcommand\rat{\mathbb Q}
\newcommand\ints{\mathbb Z}
\newcommand\comp{\mathbb C}
\newcommand\field{\mathbb F}
\newcommand\f{{\mathbb F}_q}
\newcommand\hyp{{\bf H}^2}
\newcommand\ham{\mathcal H}
\newcommand\twoorb{{\bf H}^2/\Gamma}
\newcommand\ord{\mathcal O}
\newcommand\ordtwo{\mathcal{O}'}
\newcommand\ordthree{\mathcal{O}''}
\newcommand\arith{P\rho({\mathcal O}^1)}
\newcommand\arithm{\Gamma_{\mathcal O}^1}
\newcommand\norm{\Gamma_{\ord}}
\newcommand\normone{\Gamma_{\ord}^1}
\newcommand\rf{\Gamma_{R_f}^+}
\newcommand\trace{\Gamma_{\ord}^{(2)}}
\newcommand\zmod{\mathbb Z / 2 \mathbb Z}
\newcommand\primebelow{\mathcal{p}}
\newcommand\Prime{\mathcal{P}}
\newcommand\fund{\mathcal F}
\newcommand\unit{\mathcal U}
\newcommand\units{R^{*}}
\newcommand\lattice{\mathcal L}
\newcommand\al{\alpha}
\newcommand\ga{\gamma}

\section{Introduction}
%A cocompact Fuchsian group $\Gamma$ acts properly discontinuously on $\hyp$ with quotient $\twoorb$, a hyperbolic 2-orbifold with a finite number of cone points. An arithmetic Fuchsian group has finite coarea and is therefore necessarily of the first kind.  
It is a well-known result that there are finitely many conjugacy classes of arithmetic Fuchsian groups with a given signature (\cite{MR1}, \cite{T2}).  Extensive work has been done classifying the set of $PGL_{2}(\real)$-conjugacy classes of various two-generator arithmetic Fuchsian groups: triangle groups (\cite{T1}), groups of signature (1;e) (\cite{T2}), and groups of signature $(0;2,2,2,q)$ (\cite{MR2}, \cite{ANR}).  In this paper we make progress in classifying arithmetic Fuchsian groups of signature $(2;-)$; i.e., genus two surface groups.  This is a significantly more difficult problem since these groups have much larger coarea.  

An arithmetic Fuchsian group is described by the (projectivized) group of units, $\arithm$, in a maximal order $\ord$ of a quaternion algebra over a totally real number field.  A derived arithmetic Fuchsian group is a subgroup of such a $\arithm$.  Our first main result is the classification by commensurability class of derived arithmetic Fuchsian groups of genus two, and this is summarized in Theorem \ref{list2}.  There is a finite list of signatures of groups that contain a subgroup of signature $(2;-)$.  Following \cite{MR2}, we classify all commensurability classes of derived arithmetic Fuchsian groups of the form $\arithm$ with one of these signatures by invariant quaternion algebra.  Furthermore, we determine all $PGL_{2}(\real)$-conjugacy classes of the groups $\arithm$.  Sections \ref{degree} and \ref{class} of this paper are devoted to this result and its proof.

Section \ref{max} contains our second main result, a technique for finding all $PGL_{2}(\real)$-conjugacy classes of derived arithmetic Fuchsian groups of signature $(2;-)$.  In general, the number of conjugacy classes of a subgroup of $\arithm$ is not necessarily equal to the number of $PGL_{2}(\real)$-conjugacy classes of the group $\arithm$.  However, if $\arithm$ has signature $(1;2,2)$ or $(0;2,2,2,2,2,2)$, or $(2;-)$ then index of the genus two subgroup $\Gamma$ is 4, 2, or 1, respectively.  In these cases, we can use a fundamental region along with our results from Theorem \ref{list2} to determine an explicit set of generators for $\arithm$ (using a computer program). This ultimately pins down the $PGL_{2}(\real)$-conjugacy class of the genus two subgroup $\Gamma$ as this is determined by the traces of certain products of the group generators.  Although our methods are essentially computational, we also prove some general results on the structure of maximal orders for a class of quaternion algebras associated to arithmetic Fuchsian groups.  In the last section, we use our results to explicitly determine a set of generators for a few examples of derived arithmetic Fuchsian groups of signature $(2;-)$.   
 %We remark that in \cite{T1},\cite{T2} and \cite{MR2}, the $PGL_{2}(\real)$-conjugacy classes for all of the groups $\Gamma'$ containing genus two surface groups are classified except for those of signatures 

%In Section \ref{degree}, we show that all derived genus two surface groups arise from number fields of degree less than or equal to 5.  In Section \ref{class}, we classify all commensurability classes of derived genus two surface groups by invariant quaternion algebra.  %Using Pari and estimates on the degree of the number field, we find that there exist no derived arithmetic Fuchsian groups of signature (2;-) arising from quaternion algebras over number fields of degree greater than 5.  %We list all possible quaternion algebras of the groups along with the signatures and number of conjugacy classes of their unit groups, $\arithm$.  

\section{Preliminaries}

In order to state and prove our main results, it is necessary to give a brief overview of the theory of arithmetic Fuchsian groups.  This includes a small section of number theory consisting of definitions and results that will figure prominently in our proofs.  

\subsection{Fuchsian Groups}

In this section we collect some standard results concerning Fuchsian groups.  A Fuchsian group is a discrete subgroup of $PSL_2(\real)$ that acts properly discontinuously on the hyperbolic plane $\hyp$.   Fuchsian groups of the first kind have a presentation of the form
$$
\left\langle a_1,b_1, \ldots, a_g, b_g, c_1, \ldots c_r, p_1 \ldots, p_s| \prod_{i=1}^g [a_i,b_i] \prod_{j=1}^r c_j \prod_{k=1}^s p_k , c_1^{m_1} , \ldots , c_r^{m_r}  \right\rangle,
$$
where the $c_i$ represent the $r$ conjugacy classes of maximal cyclic subgroups of order $m_i$ for $i=1, \ldots, r$.
A Fuchsian group $\Gamma$ with the above presentation has {\it signature}  
\begin{equation} \label{sig}
(g;m_1, \ldots, m_r;s).
\end{equation}
Note that $\Gamma$ is cocompact if and only if $s=0$.  Since we will be concerned only with cocompact groups, we will abbreviate the signature to $(g;m_1, \ldots, m_r)$.  A finitely generated Fuchsian group $\Gamma$ of the first kind has finite coarea, i.e., $\hyp / \Gamma$ has finite hyperbolic area, and its area can be computed using the Riemann-Hurwitz formula:
\begin{equation}
\mu(\Gamma):=\mbox{area}(\twoorb) =  2 \pi \left(2g-2 + \sum_{i=1}^r \frac{m_i-1}{m_i}+s\right). \label{rh}
\end{equation}
Furthermore, if $\Gamma_1 < \Gamma$ are Fuchsian groups and $|\Gamma: \Gamma_1| = M$ then $\mu(\Gamma_1) = M \cdot \mu(\Gamma)$.

Also, recall that two Fuchsian groups $\Gamma_1$ and $\Gamma_2$ are {\em commensurable} if they share a finite index subgroup, i.e., $|\Gamma_1: \Gamma_1 \cap \Gamma_2| < \infty$ and $|\Gamma_2: \Gamma_1 \cap \Gamma_2| < \infty$.  The {\em commensurability class} of a group $\Gamma$ is the collection of groups with which $\Gamma$ is commensurable.

\subsection{Arithmetic Fuchsian and Derived Arithmetic Fuchsian Groups}

An arithmetic Fuchsian group has finite coarea and therefore is necessarily of the first kind.   Arithmetic Fuchsian groups are defined via quaternion algebras over totally real number fields.   If $k$ is a number field and $A$ a quaternion algebra $A$ over $k$, i.e.,  a four-dimensional central simple algebra over $k$, then any quaternion algebra has an associated Hilbert symbol
$$
A=\left(  \frac{a,b}{k} \right)
$$
where $i^2 = a, j^2 = b,  ij = -ji$ for some $a, b\in k^*$.  %The basis $\{1, i,j,ij \}$ is referred to as the \emph{standard basis} of $A$.

%For any element $x = x_0 + x_1 i + x_2 j + x_3 ij \in A$, the conjugate $\bar{x}$ of $x$ is defined by $\bar{x} = x_0 - x_1 i - x_2 j - x_3 ij.$

%For $x \in A$ the (reduced) norm and (reduced) trace of $x$ lie in $F$ and are defined by 
%$$\mbox{n}(x) = x \bar{x} \quad \mbox{and}\quad \mbox{tr}(x) = x + \bar{x},$$

%respectively.  

The algebra $A$ is {\em ramified} at a real infinite place $\sigma$ of $k$ if $A \otimes_{\sigma(k)} \real \cong \ham$, where $\ham$ denotes the Hamiltonian quaternions, and {\em unramified} at $\sigma$ if $A \otimes_{\sigma(k)} \real \cong M_2(\real)$.  

Similarly, if $v$ is a finite place of $k$ and $k_v$ the completion of $k$ corresponding to $v$, then $A$ is {\em ramified} at $v$ if $A \otimes_k  k_v$ is a division algebra.  Otherwise, $A$ is {\em unramified} at $v$ and  $A \otimes_k  k_v \cong M_2(k_v)$.   

The ramification set of $A$ will be denoted by $Ram(A)$.  Furthermore, $Ram(A)=Ram_{\infty}(A)\cup Ram_f(A)$, where $Ram_f(A)$ (resp. $Ram_{\infty}(A)$) are the set of finite (resp. infinite) places at which $A$ is ramified.   We will denote the product of the primes at which $A$ is ramified by $\Delta(A)$.  

We will use the following standard results on quaternion algebras (see \cite{MRe}):
\begin{enumerate}
\item[(i)]  Let $A$ be a quaternion algebra over a number field $k$.  The number of places at which $A$ is ramified is of even cardinality.
\item[(ii)] Given a number field $k$, a collection $S_{1}=\{\sigma_{1}, \ldots, \sigma_{r}\}$ of real infinite places of $k$, and a collection $S_{2}=\{\Prime_1,\ldots, \Prime_s\}$ of finite places of $k$ such that $r+s$ is even, there exists a quaternion algebra defined over $k$ with $Ram_{\infty}(A)=S_{1}$ and $Ram_f(A)=S_{2}$.
\item[(iii)]  Let $A$ and $A'$ be quaternion algebras over a number field $k$.  Then $A \cong A'$ if and only if $Ram(A)  = Ram(A')$.
\end{enumerate}

%Furthermore, the cardinality of the full ramification set of $A$ is even, and conversely, for any set of places of $k$ of even cardinality there exists a quaternion algebra ramified at that set of places.

%Let $R$ denote the ring of integers of the number field $k$.  If V is a vector space over $k$, an {\em $R$-lattice L} is $V$ is a finitely generated $R$-module contained in $V$.  Moreover, $L$ is {\em complete} if $L \otimes_R k \equiv V$.

An {\em{order}} $\ord$ of A is a complete $R_k$-lattice which is also a ring with 1, where $R_k$ is the ring of integers in the number field $k$.  Furthermore, an order $\ord$ is maximal if it is maximal with respect to inclusion.  

%In Section 3.1, we will discuss quaternion algebras and maximal orders in more detail.

Let $k$ be a totally real field with $|k:\rat|=n$ and $A$ a quaternion algebra over $k$ which is ramified at all but one real place.   Then
$$
A \otimes_{k} \real \cong M_2(\real) \oplus \ham^{n-1}.
$$

If $\rho$ is the unique $k$-embedding of $A$ into $M_2(\real)$ and $\ord$ a maximal order in $A$, then the image under $\rho$ of the group, $\ord^1$,  of elements of norm 1 in $\ord$ is contained in $SL_2(\real)$ and the group $\arith \subset PSL_2(\real)$ forms a finite coarea Fuchsian group.  A subgroup $\Gamma$ of $PSL_2(\real)$ is an {\em arithmetic Fuchsian group} if it is commensurable with some such $\arith$.  In addition, $\Gamma$ is {\em derived from a quaternion algebra} or a {\em derived arithmetic Fuchsian group} if $\Gamma \leq \arith$.  We will denote $\arith$ by $\Gamma_{\mathcal O}^1$.   The area of $\hyp/\arithm$ can be computed by the following formula (\cite{B}):
\begin{equation}
\mbox{area}(\hyp / \arithm)             = \frac{8 \pi d_k^{3/2} \zeta_k(2) \Pi_{{\Prime}|\Delta(A)} (N({\Prime})-1)} {(4 \pi^2)^{|k:\rat|}}, \label{eq:vol}
\end{equation}
where $d_k$ is the discriminant of the number field $k$ and $\zeta_k$ is the Dedekind zeta function of the field $k$ defined for $\mbox{Re}(s) > 1$ by $\zeta_k(s) = \sum_I \frac{1}{N(I)^s}$ (the sum is over all ideals in $R_k$).

\begin{notation} Throughout the remainder of the article, we will use {\em{DAFG}} to denote a derived arithmetic Fuchsian group.
\end{notation}

If $\Gamma$ is an arithmetic Fuchsian group, then the corresponding quaternion algebra $A \Gamma$ is uniquely determined up to isomorphism and is called the {\em invariant quaternion algebra} of $\Gamma$.  Moreover, two arithmetic Fuchsian groups are commensurable if and only if their invariant quaternion algebras are isomorphic (\cite{T1}).

%Hilbert's Reciprocity Law implies the following theorem about the ramification set of a quaternion algebra over a number field:

%Since arithmetic Fuchsian groups are defined over totally real fields $k$, this implies $|Ram_{\infty}(A)|=n-1$, where $|k:\rat|=n$.  

\subsection{ Number of Conjugacy Classes}
The number of $PGL_2(\real)$-conjugacy classes of an arithmetic Fuchsian group depends on the infinite places of the number field $k$ and the number of conjugacy classes of maximal orders of the quaternion algebra $A$.  We will be solely concerned with $PGL_2(\real)$-conjugacy classes here, so  throughout the text conjugacy class should be interpreted as $PGL_2(\real)$-conjugacy class.  Most of what follows can be found in \cite{V}.

For any maximal order $\ord$ of $A$, $\norm$  will denote the arithmetic Fuchsian group
$$
\norm = \{ x \in A^* | x \ord x^{-1} = \ord \}.
$$
Let $\ord$ and $\ord'$ be two maximal orders in quaternion algebras $A / k$ and $A' / k'$, respectively.   If the groups $\norm$ and $\Gamma_{\ord'}$ are conjugate, then $k$ and $k'$ are isomorphic and 
$$
x \arithm x^{-1} = \Gamma_{\ord'}^1.
$$
A result of Vigneras (\cite{V}) states that two groups $\Gamma_{\ord}^1$ and $\Gamma_{\ord'}^1$ are conjugate if and only if there exists a $\rat$-isomorphism $\tau$ such that $\tau(A)=A'$ and $\ord' = \tau(a \ord a^{-1})$ with $a\in A$.

The class number of $k$, $h=h(k)$, is the order of the class group $I_k/P_k$, where $I_{k}$ is the group of fractional ideals of $R_k$ and $P_k$ the group of nonzero principal ideals of $R_k$.  Let $$k_{\infty}^{*}=\{x\in k | \sigma(k) > 0 \mbox{ for all } \sigma \in Ram_{\infty}(A)\}.$$  The restricted class group, whose order we will denote by $h_{\infty}$, is the group 
$$I_k/P_{k,\infty}
$$
where $P_{k,\infty}$ is the group of principal ideals with generator in $k_{\infty}^{*}$. 
We also have that 
 \begin{equation}  \label{hinf}
h_{\infty} = \frac{h 2^{n-1}}{|R_{k}^{*}:R_{k}^{*}\cap k_{\infty}^{*}|},
\end{equation}
where $R_k^{*}$ is the group of units of $R_{k}$.  The number of conjugacy classes of maximal orders in a quaternion algebra $A$ defined over $k$, denoted by $t=t(A)$, is finite and is called the {\em type number} of $A$.  It is the order of the quotient of the restricted class group of $k$ by the subgroup generated by the squares of the ideals of $R_k$ and the prime ideals dividing the discriminant $\Delta(A)$; so we have  
\begin{equation}\label{type}
t = \left| \frac{I_{k} }{ I_{k}^{2} DP_{k,\infty}} \right|,
\end{equation}
where $D$ is the subgroup of prime ideals dividing the discriminant $\Delta(A)$.  It follows that $t$ divides $h_{\infty}$.  In many cases, $h_{\infty}=1$, and we will use this to show that $t=1$.  Also, in the case that $Ram_f(A)= \emptyset$ and $h=1$, from the definitions above one can deduce that $t=h_{\infty}$.

\subsection{Torsion in Arithmetic Fuchsian Groups}

Throughout this section and the remainder of the text, $\zeta_{2m}$ will denote a primitive $2m$-th root of unity.  Also, $k_{m}$ will denote the field $\rat(\cos(\frac{\pi}{m}))$, which is the unique totally real subfield of $\rat(\zeta_{2m})$ of index 2.  Note that when $m$ is odd, $\rat(\zeta_{2m})=\rat(\zeta_m)$.  The existence of torsion in an arithmetic Fuchsian group $\arithm$ defined over a number field $k$ depends primarily on the subfields of $k$ of the form $k_{m}$ and the existence of embeddings of suitable quadratic extensions of $k$ into the invariant quaternion algebra $A$.  A more detailed treatment of this topic can be found in \cite{MRe}, Ch. 12. 

%Since the groups $\arithm$ in consideration are cocompact, the associated invariant quaterion algebra $A$ is a division algebra, and so we will always make this assumption.  
Let $A^1$ denote the elements of norm 1 in $A$ and $P(A^1)$ its projectivization.  Let  $\ord$ be a maximal order in $A$ and suppose the group $\arithm$ contains an element of order $m$.  Then $P(A^1)$ contains an element of order $m$ and $A^1$ contains an element $u$ of order $2m$.  This implies that $\mbox{tr }u \in k$ and hence $k_m \subset k$.  Furthermore, $k(\zeta_{2m})$ is a quadratic extension of $k$ that embeds in $A$.  %Conversely, if $\rat(\zeta_{2m}) | k$ embeds in $A$, then $P(A^1)$ will contain an element of order $m$.

Conversely, using the following theorem (cf. \cite{Sc}), one can show that if $k(\zeta_{2m}) | k$ is a quadratic extension that embeds in $A$, then $\arithm$ necessarily contains elements of order $m$.  %This first theorem gives necessary and sufficient conditions for the group $A^1$ to have torsion, whereas the second theorem of Chinburg and Friedman provides passage to the groups $\arithm$.  

%\begin{corollary} \label{split} Let $A$ be a quaternion algebra over a totally real number field $k$ and suppose $|k(\zeta_{2m}):k|=2$.  If there exists a prime $\Prime \in Ram_f(A)$ such that $\Prime$ splits in $k(\zeta_{2m}) | k$, then the unit group $\arithm$ of a maximal order $\ord$ will contain no elements of order $m$.
%\end{corollary}
%\begin{proof}  If $\Prime \in Ram_f(A)$ splits in $k(\zeta_{m}) | k$, then $k(\zeta_{m})  \otimes_{k} k_{\Prime}$ has zero divisors and hence $k(\zeta_{m}) $ does not embed in $A$ by the above theorem.  \end{proof}%\\

\begin{theorem} [\cite{CF}] \label{cf} Let $k$ be a number field and $A$ a quaternion division algebra over $k$ such that there is at least one infinite place of $k$ at which $A$ is unramified.  Let $\Omega$ be a commutative $R_k$-order whose field of quotients $\lattice$ is a quadratic extension of $k$ such that $\lattice \subset A$.  Then every maximal order in $A$ contains a conjugate of $\Omega$ except possibly when the following conditions both hold:

(a) $\lattice$ and $A$ are unramified at all finite places and ramified at exactly the same set of real places of $k$,

(b) all prime ideals $\Prime$ dividing the relative discriminant ideal $d_{\Omega | R_k}$ of $\Omega$ are  split in $\lattice | k$.
\end{theorem}

The order $\Omega=R_k[\zeta_{2m}]$ is a commutative $R_k$-order whose field of quotients is $\lattice=k(\zeta_{2m})$.  In the case of arithmetic Fuchsian groups, the field $k$ is totally real and the field $k(\zeta_{2m})$ is a totally imaginary extension of $\rat$.  Therefore, all real places of $k$ are ramified in $k(\zeta_{2m}) | k$; however, the algebra $A$ is ramified at all real places but one.  So condition (a) of Theorem \ref{cf} never holds.  Thus, if $\lattice \subset A$, then every maximal order $\ord$ of $A$ will contain elements of order $2m$.  Therefore, if $P(A^1)$ contains elements of order $m$, then so will $\arithm$ for any maximal order $\ord$ of $A$.  We have just shown the following:

\begin{theorem} [\cite{MRe}] \label{ramord}
An arithmetic Fuchsian group $\arithm$ contains an element of order $m$ if and only if the field $k(\zeta_{2m})$ embeds in $A$. 
\end{theorem}

The following theorem gives necessary and sufficient conditions for the embedding of the extension $k(\zeta_{2m}) | k$ into the quaternion algebra $A$.

\begin{theorem} \label{embed}
Let $A$ be a quaternion algebra over a number field  $k$ and let $ \ell | k$ be a quadratic extension.  Then $\ell$ embeds in $A$ if and only if $\ell \otimes_{k} k_v$ is a field for each $v \in Ram(A)$. 
\end{theorem} 

We can use this theorem along with Theorem \ref{cf} to give a characterization of the existence of torsion in the groups $\arithm$.  This result will be used frequently in the proof of Theorem \ref{list2}: 

\begin{lemma} \label{split} Let $k$ be a totally real number field such that $k_m \subset k$, $A$ a quaternion algebra ramified at all but one real place over $k$, and $\ord$ a maximal order in $A$. The group $\arithm$ will contain an element of order m if and only if $\Prime$ does not split in $k(\zeta_{2m})| k$ for all $\Prime \in Ram_f(A)$,

\end{lemma}
\begin{proof}  Let $\zeta_{2m}$ be a primitive $2m$-th root of unity.  By Theorem \ref{embed}, the quadratic extension $k(\zeta_{2m})$ of $k$ embeds in $A$ if and only if $k(\zeta_{2m}) \otimes_k k_v$ is a field for each $v \in Ram (A)$.  Since $k(\zeta_{2m})$ is totally imaginary, this holds automatically for all $v \in Ram_\infty(A)$.   By Theorem \ref{embed}, $k(\zeta_{2m})$ of $k$ embeds in $A$ if and only if  $\Prime$ does not split in $k(\zeta_{2m})|k$ for all $\Prime \in Ram_f(A)$.  Theorem \ref{ramord} now gives the desired conclusion. 
\end{proof}

In the case $k$ is a totally real field, the relative class number $h^{-}$ for the extension $k(\zeta_{2m}) | k$ is defined as
$$h^{-}=\frac{h(k(\zeta_{2m})) }{h(k)}\in \ints,$$ 
where $h(k(\zeta_{2m}))$ is the class number of $k(\zeta_{2m}) | \rat$ and $h(k)$ is the class number of $k$, (cf. \cite{W} p.38). 

If a maximal order $\ord$ in $A$ contains elements of finite order, then we can calculate the number of conjugacy classes, $a_m$, of maximal cyclic subgroups of order $m$ in $\arithm$, provided $\{1, \zeta_{2m}\}$ is a relative integral basis for the quadratic extension $k(\zeta_{2m})|k$ (\cite{Sc}).  If this assumption holds, then

\begin{equation}
 a_m = \frac{h^{-}}{ |R^*_{k(\zeta_{2m})}:R^{*(2)}_{k(\zeta_{2m})}|} \Pi_{{\Prime} | \Delta(A)}\left(1- \left( \frac{k(\zeta_{2m})}{{\Prime}} \right) \right),
\label{eq:tor}
\end{equation}
where $ \left( \frac{k(\zeta_{2m})}{\Prime} \right)$ is the Artin symbol (which is equal to 1, 0, or -1, according to whether $\Prime$ splits, ramifies, or is inert, respectively, in the extension $k(\zeta_{2m}) | k)$ and 
$$
R^{*(2)}_{k(\zeta_{2m})}= \{x \in R^{*}_{k(\zeta_{2m})} | N_{k(\zeta_{2m}) | k} (x) \in R_k^{*(2)}  \}.
$$
 
In some cases, we can use the following lemma to simplify the above formula (\ref{eq:tor}):

\begin{lemma}  \label{simptor} Let $k$ be a totally real number field of odd degree and suppose $\{1, \zeta_{2m} \}$ is a relative integral basis for $k(\zeta_{2m})|k$ .  If $h^{-}$ is odd, then $ |R^*_{k(\zeta_{2m})}: R^{*(2)}_{k(\zeta_{2m})}| =1$.
\end{lemma}
\begin{proof}  Both quantities $h(k(\zeta_{2m})) / h$ and $|R^{*}_{k(\zeta_{2m})}: R^{*(2)}_{k(\zeta_{2m})}|$ depend only on the number field $k$ and, hence, are independent of the quaternion algebra $A$. Since $|R^*_{k(\zeta_{2m})}: R^{*(2)}_{k(\zeta_{2m})}|$ is a finite 2-group, its order is $2^n$, for some non-negative integer $n$.  If $|k:\rat|$ is odd, then let $A$ be a quaternion algebra unramified at all finite places.  Since $a_m= \frac{h(k(\zeta_{2m}))}{h |R^*_{k(\zeta_{2m})}:R^{*(2)}_{k(\zeta_{2m})}|} \in \ints$ and $h^{-}=h(k(\zeta_{2m})) / h$ is odd, we must have $|R^*_{k(\zeta_{2m})}:R^{*(2)}_{k(\zeta_{2m})}| =1$.  
%If $|k:\rat|$ is even, then let $A$ be the quaternion algebra ramified at a prime $\Prime$ such that $\Prime$ is ramified in the extension $k(\zeta_{2m}) | k$.  Then  $\left(1- \left( \frac{k(\zeta_{2m})}{\Prime} \right) \right)=1$ and the argument follows as above. 
\end{proof}

Since $\zeta_3$ will arise frequently in our calculations, we will fix the notation $\omega = \zeta_3$.  Furthermore, the following lemma can often be used to simplify formula (\ref{eq:tor}):

\begin{lemma} \label{lem} Suppose that $k$ is a totally real number field and that 2 is unramified in $k | \rat$, then 
$$
|R^*_{k(i)}:R^{*(2)}_{k(i)}|=1. 
$$
Likewise, if 3 is unramified in $k | \rat$, then
$$
|R^*_{k(\omega)}:R^{*(2)}_{k(\omega)}|=1. 
$$
\end{lemma}

The proof of Lemma \ref{lem} requires the following two general facts (cf. \cite{Rib}, Ch. 10, and \cite{Pa}, respectively).

\begin{fact} \label{intbas} Let $k$ be a totally real number field such that $d_{k}$ is not divisible by 2.  Then $\{1,i \}$ is a relative integral basis for $k(i) |k$.  Likewise, if $d_{k}$ is not divisible by 3, then $\{1, \omega \}$ is a relative integral basis for $k(\omega) | k$. 
\end{fact}

\begin{fact}\label{units}
Let $k$ be a totally real number field and K a totally imaginary quadratic extension of $k$.  Then every unit $\epsilon$ of $K$ has the form $\epsilon = \zeta \cdot \eta$, where $\zeta$ is a root of unity with $\zeta^2 \in K$ and $\eta$ is a real unit with $\eta^2 \in k$. 
\end{fact}

\noindent
\begin{proof}[Proof of Lemma \ref{lem}]  Let us first consider the case $k(i)$.  Since 2 does not divide the discriminant of $k$, $\{1,i \}$ is a relative integral basis for the extension $k(i) | k$.  Suppose that $\cos(\frac{\pi}{m})+i \sin(\frac{\pi}{m})$ is a root of unity in $k(i)$.  Since this is also an algebraic integer, it can be written as $a + bi$, where $a, b \in R_k$.  Now, the only solutions of 
$\cos(\frac{\pi}{m})+i \sin(\frac{\pi}{m})=a + bi$ correspond to the units $\pm 1$ and $\pm i$.   By Fact \ref{units}, any unit $\epsilon$ of $k(i)$ is of the form $\epsilon = \zeta \cdot \eta$, where $\zeta^2 = \pm 1, \pm i$ and $\eta \in k$ is a real unit.  Again, using the relative integral basis, let $\epsilon = a + bi$ for some $a,b \in R_k$.   Now, any unit $\epsilon \in k(i)$ must satisfy the equation:
$$
\epsilon^2 =\zeta^2 \eta^2 =(a+bi)^2= (a^2-b^2)+2abi. 
$$
There are two possible cases to consider:\\
\noindent
\underline{Case 1}: $\pm i \eta^2= (a^2-b^2)+2abi$.\\
\noindent
\underline{Case 2}:  $\pm \eta^2  = (a^2-b^2)+2abi$.\\
\noindent
In Case 1, we must have $a = \pm b$ and $i \eta^2  =\pm 2 a^2 i$.  Since $a$ is real, this implies $\eta^2=2a^2$.  But since $a$ is an algebraic integer, $2a^2$ is not a unit.  Therefore, no unit $\epsilon$ corresponds to this case.
In Case 2, either $a=0$ or $b=0$.  Hence, $\pm \eta^2 = a^2$ or $b^2$.  Therefore, the units in $k(i)$ are of the form $\epsilon =\pm  a$ or $\epsilon = \pm bi $.  Since $\epsilon$ is a unit, this implies $a \in R_k^{*}$ or  $b \in R_k^{*}$ and, hence, $\eta^2 \in R_k^{*2}$.  Now, in either case, we have
$$
N_{k(i)|k}(\epsilon)=N_{k(i)|k}(\pm \eta)= \eta^2.
$$
This means that every unit of $k(i)$ has norm lying in $R_k^{*2}$, and so 
$$|R^*_{k(i)}:R^{*(2)}_{k(i)}|=1.$$
The proof for $k(\omega)$ is similar. \end{proof}

%Furthermore, in light of (\ref{discs}), 

It will be necessary for us to determine which periods can arise in the various number fields $k$.  First, If $k_m \subset k$, then $|k_m:\rat|$ divides $|k:\rat|$ and $d_{k_m}$ divides $d_k$.  With this in mind, we will require the following properties of the cyclotomic field $\rat(\zeta_{2p})=\rat(\zeta_p)$ and its proper subfield $k_p =\rat(\cos(\pi/p))$ when $p>3$ is a prime:

\begin{proposition} \label{cyclo}
Let $p>3$ be a prime.  Then:
\begin{enumerate}
\item $|k_p:\rat|=\frac{p-1}{2}$
\item $d_{\rat(\zeta_{p})}=p^{p-2}$ 
\item $d_{k_p} = p^{\frac{p-3}{2}}$.
\end{enumerate}
\end{proposition}
\begin{proof}  This first part follows from the fact that $|\rat(\zeta_{p}):\rat|=p-1$ and that $k_p$ is a proper subfield of index two.  The second and third parts can be found in \cite{W}, $\mbox{p.} 9 \mbox{ and p.} 28$, respectively.  \end{proof}
%These facts will be used to determine which number fields $k$ contain $k_{p}$ as a proper subfield since if $k_{p}\subset k$, then $d_{k_{p}}^{|k:k_p|}$ divides $d_k$.

\section{Bounds on the degree of the number field $k$}\label{degree}

In this section, we classify commensurability classes of DAFGs of signature $(2;-)$ by invariant quaternion algebra.  First, we determine all possible signatures of Fuchsian groups which can contain a subgroup of signature $(2;-)$.  Then, using arithmetic data, we show that all arithmetic Fuchsian groups $\arithm$ with one of these signatures are defined over a number fields of degree less than or equal to 5.

The following theorem gives necessary and sufficient conditions for the existence of torsion-free subgroups of a given index in a Fuchsian group:
\begin{theorem}[\cite{EEK}] \label{torfree}  Let $\Gamma$ be a finitely generated Fuchsian group with the standard presentation:
$$
< a_1,b_1, \ldots, a_g, b_g, c_1, \ldots c_r, p_1 \ldots, p_s| \prod_{i=1}^n [a_i,b_i] \prod_{j=1}^r c_j , \prod_{k=1}^s p_k, c_1^{m_1} , \ldots , c_r^{m_r}  >.
$$ 
Then $\Gamma$ has a torsion free subgroup of finite index $k\geq 1$ if and only if $k$ is divisible by $2^{\epsilon} \lambda$, where $\lambda = LCM(m_1, \ldots, m_r)$, and $\epsilon = 0$ if $\Gamma$ has even type, while $\epsilon = 1$ if $\Gamma$ has odd type. ($\Gamma$ has \it{odd type} if $s=0$, $\lambda$ is even, but $\lambda/m_i$ is odd for exactly an odd number of $m_i$; otherwise $\Gamma$ has \it{even type}.)
\end{theorem}

We will use this result to prove the following lemma:

\begin{lemma}  \label{list} Let $\Gamma$ be a cocompact Fuchsian group containing a genus two surface group.  Then $\Gamma$ has one of the following signatures: $(0;2,3,7)$, $(0;2,3,8)$, $(0;2,3,9)$, $(0;2,3,10)$, $(0;2,3,12)$, $(0;2,4,5)$, $(0;2,4,6)$, $(0;2,4,8)$, $(0;2,4,12)$, $(0;2,5,5)$, $(0;2,5,6)$, $(0;2,5,10)$, $(0;2,6,6)$, $(0;2,8,8)$, $(0;3,3,4)$, $(0;3,3,5)$, $(0;3,3,6)$, \\
$(0;3,3,9)$, $(0;3,4,4)$, $(0;3,6,6)$, $(0;4,4,4)$, $(0;5,5,5)$, $(0;2,2,2,3)$, $(0;2,2,2,4)$, $(0;2,2,2,6)$, $(0;2,2,3,3)$, $(0;2,2,4,4)$, $(0;3,3,3,3)$, $(0;2,2,2,2,2)$, $(0;2,2,2,2,2,2)$, $(1;2)$, $(1;3)$, or $(1;2,2)$.   
\end{lemma}

\begin{proof} If $\Gamma_1$ is a genus two surface subgroup of $\Gamma$, then $M \mu(\Gamma) = \mu(\Gamma_1) = 4 \pi$, where $M = |\Gamma:\Gamma_1|$.  This implies
$\mu(\arithm) \leq 4 \pi$.  In particular, the genus $g$ of $\Gamma$ must be less than or equal to $2$.    
Furthermore, by Theorem \ref{torfree},  depending on the signature of the group $\Gamma$, either $\lambda$ or $2\lambda$ divides the index $M$, where $\lambda = LCM(m_1, \ldots, m_r)$.  This gives us bounds on the possible torsion of $\Gamma$.  In particular, for fixed $g$, this gives us an upper bound on the number of conjugacy classes of elliptic elements:
\begin{enumerate}
\item[(1)] If $g=0$, then $\Gamma$ has at most six conjugacy classes of elliptic elements.
\item[(2)] If $g=1$, then $\Gamma$ has at most two conjugacy classes of elliptic elements.
\item[(3)] If $g=2$, then $\Gamma$ has no elliptic elements and $\Gamma=\Gamma_1$.
\end{enumerate}

For example, suppose $g=0$ and $\Gamma$ has four conjugacy classes of elliptic elements $x_i$ of order $m_i$, $1 \leq i \leq 4$.  By the Riemann-Hurwitz formula (\ref{rh}),
$$
\mu(\Gamma) = 2\pi \left(-2+\sum_{i=1}^4\frac{m_i-1}{m_i}\right).
$$
Therefore, if $\Gamma$ contains a torsion-free subgroup of genus two, then
$$
M \mu(\Gamma) = 2M\pi \left(-2+\sum_{i=1}^4\frac{m_i-1}{m_i} \right)= 4 \pi = \mu(\Gamma_1). 
$$
This translates to the existence of integers $M$, $m_i$, $1 \leq i \leq 4$ satisfying the equation:
\begin{equation}\label{sub} 
\sum_{i=1}^4\frac{m_i-1}{m_i}=\frac{2}{M}+2.
\end{equation}
In addition, $\lambda = LCM(m_1, \ldots, m_4)$ divides $M$.  Since $\mu(\Gamma) >0$, there exists at least one $x_i$ with order $m_i > 2$.  Also, we can deduce the following two facts.
\begin{enumerate}
\item[(i)] There cannot exist more than two distinct $m_i$ corresponding to the $x_i$.
\item [(ii)] If $m_i > 2$ for all $1 \leq i \leq 4$, then $m_1 = \cdots =m_4=3$.   
\end{enumerate}
If (i) does not hold, then $m_1 \geq 2$, $m_2 \geq 3$, $m_3 \geq 4$, and $M \geq \lambda \geq 6$, and this gives the following contradiction:
$$
 \frac{29}{12} \leq \sum_{i=1}^4\frac{m_i-1}{m_i}=\frac{2}{M}+2 \leq   \frac{28}{12}.
$$
Similarly, if (ii) does not hold, then $M \geq \lambda \geq 4$, and we arrive at the contradiction
$$
\frac{10}{4} \leq \sum_{i=1}^4\frac{m_i-1}{m_i}=\frac{2}{M}+2 \leq   \frac{9}{4}.
$$
Without loss of generality, suppose $x_1 \leq  x_2 \leq \cdots \leq x_4$.

\noindent
\underline{Case 1}:  $m_1=m_2=m_3=2$ and $m_4=m > 2$.

\noindent
In this case, equation (\ref{sub}) becomes
$$
\frac{3}{2}+\frac{m-1}{m}=\frac{2}{M}+2 \Longleftrightarrow \frac{m-1}{m} = \frac{2}{M}+\frac{1}{2}.
$$
The solution $m=M=6$ gives the maximal value of $m$.  The only other solutions in this case occur when $(m,M)=(3,12),(4,8)$.  

\noindent
\underline{Case 2}: $m_1=m_2=2$ and $m=m_3=m_4 > 2$.

\noindent
Again, equation (\ref{sub} becomes
$$
1+\frac{2(m-1)}{m}=\frac{2}{M}+2 \Longleftrightarrow \frac{2(m-1)}{m} = \frac{2}{M}+1.
$$
\noindent
The case $m=N$ gives the maximal value for $m$ and this occurs when $m=M=4$.  The only other possible solution occurs when $(m,M)=(3,6)$.

\noindent
\underline{Case 3}: $m_1=2$ and $m=m_2=m_3=m_4 > 2$.

\noindent
Equation (\ref{sub}) translates to
$$
\frac{1}{2}+\frac{3(m-1)}{m}=\frac{2}{M}+2;
$$
and one can easily verify that there exist no integer solutions to this equation.

\noindent
\underline{Case 4}: $m=m_1=m_2=m_3=m_4 > 2$.

\noindent
In this situation,
$$
\frac{4(m-1)}{m}=\frac{2}{M}+2,
$$
and $m=M=3$ is the only solution.  
The existence of a torsion-free subgroup of index $M$ for a group of fixed signature is guaranteed by Theorem \ref{torfree}.   Therefore, the only Fuchsian groups with signature $(0;x_1,x_2,x_3,x_4)$ containing a torsion-free subgroup of genus 2 are those with signatures $(0;2,2,2,3)$, $(0;2,2,2,4)$, $(0;2,2,2,6)$, $(0;2,2,3,3)$, $(0;2,2,4,4)$, and $(0;3,3,3,3)$.
In this manner, we analyze torsion in groups of a fixed signature to obtain the list in the Lemma.
\end{proof}

\begin{proposition}  There exist no DAFGs of signature $(2;-)$ arising from quaternion algebras over number fields of degree greater than 5.
\end{proposition}

\begin{proof} If $\arithm$ contains a genus two surface group $\Gamma$ of index $M=|\arithm:\Gamma|$, then 
\begin{equation}\label{index}
\mu(\Gamma)=4 \pi =M \mu(\arithm)=M\frac{8 \pi d_k^{3/2}\zeta_k(2) \Pi_{\Prime|\Delta(A)} (N(\Prime)-1)} {(4 \pi^2)^{|k:\rat|}}. 
\end{equation}
In particular, this implies
\begin{equation} \label{ineq}
4 \pi \geq \frac{8 \pi d_k^{3/2}\zeta_k(2) \Pi_{\Prime|\Delta(A)} (N(\Prime)-1)} {(4 \pi^2)^{|k:\rat|}}. 
\end{equation}

Note that 
\begin{equation}\label{even}
\zeta_k(2) \prod_{\Prime|\Delta(A)}(N(\Prime)-1)> \prod_{\Prime|\Delta(A)} \frac{(N(\Prime))^2}{(N(\Prime)+1)} \geq \left\{ \begin{array}{ll}
1 & \mbox{if $|k:\rat|$ is odd} \\
4/3 & \mbox{if $|k:\rat|$ is even}
\end{array} \right..
\end{equation}

Using $\Pi_{\Prime|\Delta(A)}(N(\Prime)-1)\geq 1$ and $\zeta_k(2) \geq 1$ in the above inequality gives 
\begin{equation}\label{1stest} 
4 \pi \geq \frac{8\pi d_k^{3/2}}{(4\pi^2)^{|k:\rat|}}.
\end{equation}
We now use Odlyzko's lower bounds (cf. \cite{O}) on the discriminant of a totally real number field to get an upper bound on the degree of $k$:
$$
|d_k| \geq (2.439 \times 10^{-4})(29.099)^{n},
$$
where $n=|k:\rat|$.  Using these estimates in inequality (\ref{1stest}) gives $n \leq 8$.
However, the smallest discriminants of a totally real field of degree 7 and 8 are 20,134,393 and 282,300,416, respectively (\cite{Co} and \cite{octic}).  In each case, inequality (\ref{ineq}) is violated:
$$
4 \pi \geq \frac{8\pi d_k^{3/2}}{(4\pi^2)^{|k:\rat|}} \geq \left\{ 
\begin{array}{l}
\frac{8 \pi (20,134,393)^{3/2}}{(4\pi^2)^7} \simeq 15.1925  > 4\pi \\
\frac{8 \pi (282,300,416)^{3/2}}{(4\pi^2)^8} \simeq 20.2036  > 4 \pi
\end{array} \right..
$$
\noindent
Hence, there cannot exist a DAFG of signature $(2;-)$ if $|k:\rat| \geq 7$.

To eliminate the case $|k: \rat| =6$, we again exploit the area formula (\ref{eq:vol}) and inequality (\ref{even}) to get the following inequality:
$$
\mu(\Gamma) =4 \pi \geq |\arithm:\Gamma| \frac{8 \pi d_k^{3/2} \zeta_k(2)\Pi_{{\Prime}|\Delta(A)} (N({\Prime})-1)} {(4 \pi^2)^6}
\geq \frac{ 32 \pi d_k^{3/2} } {3 (4 \pi^2)^6}.
$$
This gives us the following upper bound on the discriminant $d_k$:
\begin{equation}\label{deg6}
d_k \leq \left(\frac{3(4 \pi^2)^6}{8}\right)^{2/3} < 1,263,165.
\end{equation}
According to the lists from \cite{Co}, there are 20 number fields $k$ of degree 6 satisfying the above inequality.  For each field $k$, we investigate the behavior of small primes and, if necessary, estimate $\zeta_k(2)$ using Pari.  Since $n=6$, $|Ram_f(A)| \neq \emptyset$.  Therefore, $\Pi_{\Prime|\Delta(A)} (N(\Prime)-1) \geq N(\Prime_0)-1$, where $\Prime_0$ is the prime of smallest of norm in $k$. 

For example, consider the totally real field $k$ of degree 6 and discriminant $d_{k}~=~722,000$.  A minimal polynomial for $k$ is $f(x)= x^6-x^5-6x^4+7x^3+4x^2-5x+1$.  By Pari, the prime of smallest norm in $R_{k}$ is the unique prime $\Prime$ lying over 2 with $N(\Prime)=4$.  This implies that any group $\arithm$ defined over $k$ has area at least
$$
\frac{8 \pi d_k^{3/2} \zeta_k(2)\cdot 3} {(4 \pi^2)^6}=\frac{21\pi}{5} > 4\pi.
$$
Hence, there exist no DAFGs of signature $(2;-)$ defined over $k$.   In this fashion, we obtain a contradiction to the inequality $\mu(\Gamma) \leq 4 \pi$ for each totally real field $k$ of degree 6 with discriminant $d_k$ satisfying (\ref{deg6}). \end{proof}

\noindent
\section{Classification by Commensurability Class}\label{class}

In this section, we classify DAFGs of signature $(2;-)$ by invariant quaternion algebra.  All the groups in Lemma \ref{list} except those with one of the three signatures $(1;2,2)$, $(0;2,2,2,2,2,2)$, and $(2;-)$, have commensurability classes that have already been classified; i.e., they are  all commensurable with an arithmetic Fuchsian triangle group or one having signature $(1;e)$ or $(0;2,2,2,e).$  So it suffices to classify the commensurability classes of the groups $\arithm$ of these remaining three types and to extract the relevant results from \cite{ANR}, \cite{MR2}, \cite{T1}, and \cite{T2}.

The proof classification is exhaustive.  For each fixed degree $|k:\rat|$, we use equation (\ref{eq:vol}) to get upper bounds on the discriminant of the number field $k$.   Then we determine the existence of all quaternion algebras whose unit groups have one of the above three signatures.  Rather than go through an analysis of each number field that can correspond to such an arithmetic Fuchsian group, we give an idea of the overall approach by a few illustrative examples.  Our argument will be organized by the degree of the number field.

We will make extensive use of the following lemma, which is particularly useful in the case $|k:\rat|$ odd (since we can have $Ram_f(A)=\emptyset$ in this case):

\begin{lemma}  \label{sixdivides} If $A$ is a quaternion algebra defined over a totally real field ramified at all but one infinite place and unramified at all finite places, then $\arithm$ contains elements of orders 2 and 3.  Furthermore, if $\Gamma$ is a genus two surface group contained in $\arithm$ for $\ord$ a maximal order in $A$, then $6$ divides $|\arithm: \Gamma|$. 
\end{lemma}
\begin{proof}  Since $ Ram_f(A)$ is empty, by Lemma \ref{split}, there is no obstruction to embedding $\mathcal{L}$  in $\ord$, where $\mathcal{L} \cong {\rat}(\mbox{i})$ or $\rat(\omega)$.  Therefore, any order $\ord$ in $A$ will contain elements of orders 2 and 3.  By Theorem \ref{torfree}, if $\arithm$ has signature $(g;x_1,\ldots, x_r)$ and $\Gamma \subset \arithm$ is torsion-free, then  $6$ divides  $\lambda$ which in turn divides $|\arithm: \Gamma|$. \end{proof}

\subsection{Quintic Number Fields}

\begin{lemma}
For $|k:\rat|=5$, the only DAFGs of signature $(2;-)$ arise from quaternion algebras over the totally real fields of discriminants $d_k = 38569$, $36497,$ and $24217$.  
\end{lemma}

\begin{proof} Suppose there exists a DAFG $\Gamma < \arithm$ of genus two which is torsion-free and defined over a totally real quintic number field $k$. Using the inequalities $\zeta_k(2) \geq 1$, $\prod_{\Prime| \Delta(A)}(N(\Prime)-1) \geq 1$ in conjunction with (\ref{eq:vol}) we get that $$d_k \leq 131,981.$$ 

However, if $Ram_f(A)= \emptyset$, the index $M=[\arithm:\Gamma] \geq 6$ by Lemma \ref{sixdivides}.  Substituting back into the area formula (\ref{eq:vol}) gives $d_k \leq 39970$.  For those fields with  $39970  \leq d_k \leq 131981$, we analyze the behavior of small primes in $k$ to determine the possible ramification sets for each field $k$.  According to \cite{Co}, there are 15 number fields with $d_k < 131981$. In a few cases, small primes do exist, but we can eliminate these cases using torsion. 

For example, let $k$ be the number field with discriminant $d_k = 106069$.  Using the minimal polynomial $f(x)=x^5-2x^4-4x^3+7x^2+3x-4$ to generate the $k$ in Pari, we compute that
$$
\frac{8 \pi \cdot  106069^{3/2} \zeta_k(2)}{(4 \pi^2)^5} = 4 \pi.
$$
Furthermore, there exists a unique prime $\Prime_2$ of norm 2 in $R_k$.  Together with the fact that  $|Ram_f (A)|$ is even, this implies that $Ram_f(A) = \emptyset$.  So, by (\ref{eq:vol}), $\mu(\arithm) = 4 \pi$ for any maximal order $\ord$ in $A$.  However, by Lemma \ref{sixdivides}, $\arithm$ contains elements of orders 2 and 3.  Hence, $\arithm$ is not torsion-free, and since $\mu(\arithm) = 4 \pi$, it neither is a genus two surface group nor does it contain a genus two subgroup. 

The case $d_k = 38569$ yields a positive result.  By Pari, using the minimal polynomial $f(x)= x^5-5x^3+4x-1$, we compute that
$$
\mu(\arithm)=\frac{8 \pi  \cdot 38569^{3/2} \zeta_k(2)}{(4 \pi^2)^5} = \frac{2 \pi}{3}.
$$
So if $\Gamma \subset \arithm$ has signature $(2;-)$, then the following equation must be satisfied:
\begin{equation}
 4 \pi = \mu(\Gamma) = M \mu(\arithm) =M \frac{ 2 \pi \prod_{\Prime | \Delta(A)} (N(\Prime)-1) }{3}, 
\end{equation}
where $M=|\arithm: \Gamma|.$

Again the only solution to the above equation occurs when $M=6$ and $Ram_f(A) =\emptyset$.  Since $d_k= 38569$ is prime, $k$ contains no proper subfields other than $\rat$.  Thus, the only possibilities for elements of finite order in $\ord$ are 2 and 3.  By Lemma \ref{sixdivides}, any group $\arithm$ contains elements of orders 2 and 3.  So, in this case we get
$$
\mu(\arithm) =\frac{2 \pi}{3} = 2 \pi \left( 2g-2 +\frac{a_2}{2} + \frac{2 a_3}{3} \right). 
$$
We see that the only solution to this equation is $a_2 = a_3 =2$.  Hence, $\arithm$ has signature $(0;2,2,3,3)$ in this case, and Theorem \ref{torfree} guarantees the existence of a torsion-free subgroup of index 6.  

The totally real number fields of degree 5 with discriminants 36497 and 24217 are the only other fields which yield positive existence results. \end{proof}

\subsection{Quartic Number Fields}
%For $[k : \rat] = 3,4$, there are more possibilities for torsion.  This complicates matters slightly, but in many cases, we can still use Theorem \ref{eq:tor} to determine the signature of $\arithm$. 

For $|k: \rat|=4$, $Ram_f(A) \neq \emptyset$.  Using Proposition \ref{cyclo} to analyze the cyclotomic extensions with degree dividing 8, one can easily show that 2, 3, 4, 5, 6, 8, 10, 12, 15 are the only possible cycles for elliptic elements in this case.  %However, there are more possibilities for torsion as we now show.

%\begin{lemma}  \label{quartor} Let $k$ be a totally real number field with $|k:\rat|=4$ and denote by $k_m$ the field $\rat(\cos(\frac{2\pi}{m}))$.  Then  $|k_{m}:\rat|$ has degree dividing 4 for the following values of m: 2,3,4,5,6,8,10,12,15.  
%\end{lemma}
%\begin{proof}  The periods 2 and 3 are always possible since $\rat$ is a proper subfield of $k$.  The quadratic fields $k_4=\rat(\sqrt{2})$, $k_5=\rat(\sqrt{5})$, $k_6=\rat(\sqrt{3})$ can be possible proper subfields of $k$, which yield $4,5,6$.   We examine the cyclotomic fields $\zeta_{2m}$ denote a primitive $m$-th root of unity for which $|\rat(\zeta_{2m}):\rat|=4,8$.  If $m$ is odd, then $k_{m}$ will be a proper subfield of index two, $|k_{m}:\rat|=4$.  If $m$ is even, we need to account for the fact that $-\mbox{Id} \sim \mbox{Id}$ in $PSL_2(\real)$; this means $|k_m:\rat|=\mbox{deg}\Phi_m$, where $\Phi_m(x)$ is the $m$-th cyclotomic polynomial.  By Prop. \ref{cyclo}, there are no rational primes $p$ with $|\rat(\zeta_p):\rat|=4$.  If $m$ is a prime power $p^k$, then $\mbox{deg}\Phi = p^{k-1}(p-1)=4$. The only solution in this case occurs when $p=2$ and $k=3$ and, hence, 8 is a possible period.  If $m=p_1^{k_1}\cdots p_r^{k_1}$, then $\mbox{deg}\Phi_m= p_1^{k_1-1}(p_1-1)\cdots p_r^{k_r-1}(p_r-1)$.  Analyzing solutions to this equation when $\mbox{deg} \Phi_m=4,8$ yield the periods 10,12,15.   This completes the proof of the lemma.
%\end{proof}

Using equation (\ref{index}) in conjunction with inequality (\ref{even}), we obtain the following inequality when $|k:\rat|=4$: 
$$
4 \pi \geq \frac{32 \pi d_k^{3/2}}{3(4\pi^2)^4}.
$$
Therefore,
$$
d_k \leq \left( \frac{3(4\pi^2)^4}{8}\right)^{2/3} <  9397.
$$
There are 48 number fields with discriminants satisfying the above inequality in \cite{Co}.  %We list these number fields in Appendix A.   
Again, we eliminate all fields except those listed in Theorem \ref{list2} by estimating $\zeta_k(2)$  and examining the factorization of small primes using Pari.  

\begin{lemma}  There exist no DAFGs of signature $(2;-)$ defined over the totally real field $k$ with $d_k=5744$. 
\end{lemma}
\begin{proof}  The minimal polynomial for $k$ is $f(x)=x^4-5x^2-2x+1.$  Using Pari, we compute that
$$\mu(\arithm)=\frac{8 \pi \cdot 5744^{3/2 }\zeta_k(2)}{(4 \pi^2)^4} =\frac{5\pi}{3}.$$  
Hence, a torsion-free genus two subgroup $\Gamma$ of index $M=|\arithm:\Gamma|$ corresponds to a solution of the equation
$$
 \frac{5M}{3} \prod_{\Prime| \Delta(A)} (N(\Prime)-1) = 4.
$$
But as $M,N(\Prime) \in \ints$, this clearly has no solution. \end{proof}

We now list some positive results for the case $|k:\rat|=4$.

\begin{lemma}
Let $k$ be the totally real number field with $d_k = 3981$.  Then the only DAFG of signature $(2;-)$ defined over $k$ has invariant quaternion algebra $A$ with $Ram_{f}(A)=\Prime_{3}$, where $\Prime_{3}$ is the unique prime in $R_k$ lying over 3.  Furthermore, there is only one conjugacy class of DAFGs of this signature defined over $k$.
\end{lemma}
\begin{proof} The number field $k$ is equal to $\rat(\alpha)$ where $\alpha$ is a root of the polynomial $f(x)=x^4-x^3-4x^2+2x+1$.  Since $d_k=3 \cdot 1327$, $k$ contains no other proper subfield other than $\rat$;  the only possible non-trivial elements of finite order of $\arithm$ are those of order 2 or 3.  By Pari, we compute
$$ \mu(\arithm)=\frac{8 \pi  \cdot 3981^{3/2 }\zeta_k(2)}{(4 \pi^2)^4} =\pi.$$  
Therefore, if $\arithm$ contains a subgroup $\Gamma$ of signature $(2;-)$ of index $M=|\arithm:\Gamma|$, then
\begin{equation}\label{no4}
 M \prod_{\Prime| \Delta(A)} (N(\Prime)-1) = 4.
\end{equation}
This implies that $N(\Prime) \leq 5$ for any prime $\Prime \in Ram_f(A)$.  By Pari, we find that there are only two primes in $R_k$ with norm less than 5: $\Prime_3$ and $\Prime_5$ with $N(\Prime_3)=3$ and $N(\Prime_5)=5$.   

Thereofre $M=1,  Ram_f(A) = \{ \Prime_5 \}$ is a possible solution to (\ref{no4}).  However, $\Prime_5$ is inert in $k(\omega) | k$ which, by Lemma \ref{split}, implies that $a_3 \neq 0$.  However, by Theorem \ref{torfree}, $3 | M=1$, which is a contradiction.

We also have $M=2$ and $Ram_f(A) = \{ \Prime_3 \}$ as a possible solution to (\ref{no4}).  Furthermore, $ 6$ does not divide $d_k$, so by Lemma \ref{intbas} we can calculate $a_2$ and $a_3$ using (\ref{eq:tor}).   Using Pari, we compute that $\Prime_3$ is inert in the extension $k(i)|k$ that $h(k(i)) = 3$ and $h(k)=1$.  Also, since $2$ does not divide $d_k$, by Lemma \ref{lem} we have that $|R^*_{k(i)}:R^{*(2)}_{k(i)}|=1.$  Therefore
$$
a_2 = \frac{h^{-}}{ |R^*_{k(i)}:R^{*(2)}_{k(i)}|} \Pi_{{\Prime} | \Delta(A)}\left(1- \left( \frac{k(\zeta_{i})}{{\Prime}} \right) \right) =  3 \cdot 2 = 6.
$$
The prime $\Prime_3$ splits in $k(i)$, so by either Lemma \ref{split} or equation (\ref{eq:tor}), $a_3 = 0$.

Since these are the only possible periods for this number field, $\arithm$ must have signature $(0;2,2,2,2,2,2)$.  Finally, $\arithm$ contains a torsion free subgroup $\Gamma$ of index 2 by Theorem \ref{torfree}; since $\mbox{vol}(\Gamma)=4\pi$, $\Gamma$ must have signature $(2;-)$.  

Since the extension $k | \rat$ is not Galois and $k$ contains no proper subfields other than $\rat$, the groups corresponding to the various infinite places of $k$ will each contribute at least one conjugacy class.  For each of these quaternion algebras, we determine the type number by analyzing the embeddings of the units.  By Pari, a fundamental system of $R_k^*$ is $\{-1,\al,\al-1,\alpha^2+ \al -1\}$.  The signs of these generators at the various embeddings are shown in the table below:
$$
\begin{array}{c|ccc}
& \alpha & \alpha -1 & \alpha^2+ \al -1 \\ \hline
\alpha_{1}\simeq-1.7508 & -& - &  +\\
\alpha_{2}\simeq-0.3184 & -&  -& - \\
\alpha_{3}\simeq0.7853& + &  -& +\\ 
\alpha_{4}\simeq2.2840 &+ &  +&  +\\
\end{array}
$$
For each choice of unramified real place $\sigma_i$, $h_{\infty}=1$.  Hence, there are four distinct conjugacy classes of groups of signature $(0;2,2,2,2,2,2)$ defined over $k$.
\end{proof}

\begin{lemma}  Let $k=\rat(\sqrt{2},\sqrt{3})$ be the number field with $d_{k}=2304$.  The only DAFG of signature $(2;-)$ arising from a quaternion algebra $A$ over $k$ has $Ram_{f}(A) = \Prime_{3}$ where $\Prime_{3}$ is the unique prime of norm 9 in $k$.
\end{lemma}
\begin{proof}  The periods 2, 3, 4, and 6 are all obvious possibilities for torsion since $\rat$, $k_2$, $k_3$ are proper subfields of $k$.  The fact that $5 \not| 2304 = 2^8 \cdot 3^2$ implies that these are the only possibilities.  In this case, the $k=\rat(\al)$, where $\al$ is a root of the polynomial $f(x)=x^4-4x^2+1$.  Using Pari, we compute that
$$\frac{8\pi \cdot 2304^{3/2}\zeta_k(2)}{(4 \pi^2)^4} =\frac{\pi}{2}.$$  
Therefore, the existence of a torsion-free genus two subgroup amounts to the existence of a solution to the equation
\begin{equation}\label{no5}
 M \prod_{\Prime| \Delta(A)} (N(\Prime)-1) = 8.
\end{equation}
The only primes $\Prime$ in $R_k$ with $(N(\Prime)-1)$ dividing 8 are the unique primes $\Prime_2$ and $\Prime_3$ of norms 2 and 9, respectively.  Since $|Ram_f(A)|$ is odd, the only solution to (\ref{no5})
is $M=1$ and $Ram_f(A)=\{ \Prime_3\}$.  The prime $\Prime_3$ splits in both $k(i) | k$ and $k(\omega) | k$.  So, by Lemma \ref{split}, for any maximal order $\ord$, the group $\arithm$ contains no elements of orders 2 or 3.  This also implies that $\arithm$ contains no elements of order 4 or 6; therefore, $\arithm$ is torsion-free and has genus 2.  Since $k | \rat$ is Galois, there is only one conjugacy class of arithmetic Fuchsian groups $\arithm$ defined over $k$.\end{proof}

\begin{lemma}  Let $k$ be the number field with $d_{k}=1957$.  Then the only DAFGs of signature $(2;-)$ arising from a quaternion algebra $A$ over $k$ containing genus two subgroups are those listed in Theorem \ref{list2}.
\end{lemma}

%\begin{proof}  
The number field $k$ with discriminant $d_{k}$ is equal to $\rat(\alpha)$ where $\alpha$ is a root of the polynomial $f(x)=x^{4}-4x^{2}-x+1$. Since $d_k=1957=19\cdot 103$, $k$ contains no proper subfields other than $\rat$.  Using Pari, we compute that
$$
\frac{8\pi \cdot 1957^{3/2}\zeta_{k}(2)}{(4\pi^{2)^{4}}}=\frac{\pi}{3}.
$$
Again, we consider solutions to the equation
$$4\pi=M \frac{\pi  \prod_{\Prime | \Delta(A)} (N(\Prime)-1)}{3},$$
or equivalently,
\begin{equation} \label{1957sol} 
M \prod_{\Prime | \Delta(A)} (N(\Prime)-1) = 12, 
\end{equation}
by analyzing the primes in $R_{k}$.  In particular, any prime $\Prime$ in the ramification set of $A$ has norm at most 13.  The rational primes 2, 5, and 13 remain prime in the extension $k| \rat$, so they cannot lie in $Ram_{f}(A)$.  By Pari, there are two primes $\Prime_{3}=(\al-2)R_k$ and $\Prime_{3}'$ lying over 3, with norms $N(\Prime_{3})=3$ and $N(\Prime_{3}')=9$, respectively.  There also exists a prime $\Prime_{7}=(2 \alpha +1)R_{k}$ of norm 7.   Combining this with the fact that $|Ram_{f}(A)|$ is odd, we get that the only possible solutions $(M,Ram_{f}(A))$ to the above equation are $(6,\Prime_{3})$ and $(2,\Prime_{7})$.  The quaternion algebras with $Ram_{f}(A)=\{\Prime_{3}\}$ appear in the lists of \cite{MR2} and the unit groups $\arithm$ are of signature $(0;2,2,2,3)$ in this case.  

Let us consider the algebra with $Ram_{f}(A)=\{\Prime_{7} \}$. Since $6\not| 1957$ and $k$ contains no proper subfields other than $\rat$, $\arithm$ can only have elements of orders 2 and 3.  Moreover, we can compute the number of elements of orders 2 and 3 using formula ({\ref{eq:tor}).  Since $\Prime_{7}$ splits in $k(\omega)|k$, $\arithm$ contains no elements of order 3.  Since $\Prime_{7}$ is inert in $k(i) | k$, and $h(k(i))=1$, there are two conjugacy classes of elements of order 2.  Therefore, any group $\arithm$ arising from $A$ has signature $(1;2,2)$.

In order to determine the number of conjugacy classes of the groups $\arithm$, we again analyze the behavior of the units $R_{k}^{*}$ at the various embeddings $\alpha_{i}$.  The set $\{-1, \alpha,\alpha-1,\alpha+ 2 \}$ is a fundamental system of units for $R_{k}^{*}$ and the following table lists the signs of the generators:
$$
\begin{array}{c|ccc}
& \alpha & \alpha -1 & \alpha +2 \\\hline
\alpha_{1}\simeq -2.0615 & -& - & - \\
\alpha_{2}\simeq -0.3963 & - & - &+\\
\alpha_{3}\simeq 0.6938 &+ & - & +\\ 
\alpha_{4}\simeq 1.7640 &+& + & + \\
\end{array}.
$$
Since the extension $k | \rat$ is not Galois and $k$ contains no proper subfields, we again get at least one conjugacy class corresponding to the algebra unramified at the place $\alpha_{i}, 1 \leq i \leq 4$.  The class number of $k$ is 1, so $h_{\infty}=2^3/|R_{k}^{*} / R_{k}^{*} \cap k_{\infty}^{*}|$.  By the table above, we see that for each choice of $\sigma_i$,  $h_{\infty}=1$ for the algebra unramified at $\sigma_i$.  Therefore, there are exactly four conjugacy classes of groups of signature $(1;2,2)$ arising from quaternion algebras defined over $k$. 
%\end{proof}

%\begin{lemma}  Let $k$ be the number field with $d_{k}=2000$.  Then the only derived group arising from a quaternion algebra $A$ over $k$ satisfies $Ram_{f}(A) =\Prime_{2}$, where $\Prime_{2}$ is the unique prime of $k$ dividing 2.
%\end{lemma}
%\begin{proof} Since $d_k=2000$, $\lambda_5 \subset k$.  However, we cannot use formula (\ref{eq:tor}) to determine the periods of $\arithm$, but this is not too problematic for us.
%$$
%m \frac{\pi  \prod_{\Prime | \Delta(A)} (N(\Prime)-1)}{3} = 4 \pi,
%$$
%or, equivalently, 
%\begin{equation}
%m \prod_{\Prime | \Delta(A)} (N(\Prime)-1) = 12 \label{eq:no2}
%\end{equation}
%must hold.  A possible solution to (\ref{eq:no2}) is $m=3, Ram_f(A) = \{ \Prime_5 \}$. Certainly, $3 \mid d_k$, so $\{1, \omega \}$ is a relative integral basis for $k(\omega) / k$, and $a_3 = 4$ as $h(k(\omega))=2, h=1$, and  $\Prime_5$ is inert in $k(\omega / k)$. Now, $\{1,i \}$ and $\{1,\lambda_5 \}$ are not relative integral bases for $k(i) /k $ and $k(\lambda_5) /k $, respectively, but one easily sees that the only solution to 
%$$
%\pi = 4 \frac{4\pi}{3} = 2 \pi \left( 2g - 2 + \frac{a_2}{2} + 4 \frac{2}{3} + \frac{4 a_5}{5} \right)
%$$
%is $a_2=a_5=0, g=0$.

%Since there are no primes of norm $\leq 4$, the only other possible solution to (\ref{eq:no2}) is $m=4,  Ram_f(A) = \{ \Prime_2 \}$, where  $\Prime_2$ is the unique prime of norm 4 in $k$.  However, $\Prime_2$ does not split in $k(\lambda_5) / k$.  Therefore, by Theorem \ref{cf}, $a_5 \neq 0$.  Hence, $5 | m=4$, which is a contradiction.  \end{proof}
\subsection{Cubic Number Fields}

\begin{lemma}  The only possibly periods of elements of finite order that can arise in $\arithm$ defined over fields $k$ with $|k: \rat|=3$ are 2,3,7, and 9.
\end{lemma}
\begin{proof}  Since $k$ contains no proper subfields other than $\rat$, 2 and 3 are the only possible periods than can arise from proper subfields of $k$.  By Prop. \ref{cyclo}, 7 is the only prime for which $|k_{p}:\rat|=3$.  In fact, $k_7=\rat(\cos(\frac{\pi}{7}))$ is the totally real cubic field with discriminant 49.  For prime powers $m=p^k$, the only field $\rat(\zeta_{2m})$ with $|\rat(\zeta_{2m}):\rat|=6$ is $m=9$.  This corresponds to the totally real field of discriminant 81.  There are no composite $m$ for which $|\rat(\zeta_{2m}):\rat|=6$; this finishes the proof. 
\end{proof} 
If $[k : \rat] =3$, it is possible that $Ram_f(A) = \emptyset$.  As in the case $[k : \rat] = 5$, this helps to simplify the process immensely, since this implies $d_k \leq 297$.  If $Ram_f(A) \neq \emptyset$, then $d_k \leq 981$.   On the lists \cite{Co}, there are 25 number fields $k$ with discriminants satisfying the latter inequality (see \cite{M} Appendix A).  The cases $k=\rat(\cos(\frac{\pi}{9}))$ and $k=\rat(\cos(\frac{\pi}{7}))$ where 9 and 7, respectively, are possible periods require special examination.  We analyze the latter case in detail below.

\begin{lemma}  There exist no DAFGs of signature $(2;-)$ arising from a quaternion algebra defined over the totally real cubic number field of discriminant 361.
\end{lemma}

\noindent
\begin{proof}  The field $k$ has minimal polynomial $f(x)=x^3-x^2-6x+7$.  Using Pari, we compute 
$$
\mu(\arithm)=\frac{8 \pi \cdot 361^{3/2}\zeta_k(2)}{(4 \pi^2)^3} = \pi.
$$
By the preceding comments, $Ram_f(A)\neq \emptyset$.   Now the equation
\begin{equation}
 4 \pi = \mu(\Gamma) =M \mu(\arithm) =M \pi \prod_{\Prime | \Delta(A)}( N(\Prime)-1) 
\end{equation}
has no solutions, as 2, 3, and 5 are inert in $k | \rat$.  Therefore, there are no DAFGs 
of signature $(2;-)$ defined over $k$.
\end{proof} 
\noindent
\begin{lemma} For $k=\rat(\cos(\frac{\pi}{7}))$, the only possible $\arithm$ containing a subgroup $\Gamma$ of signature $(2;-)$ are those listed in Theorem \ref{list2}.  
\end{lemma}
\noindent
\begin{proof} Fix $f(x)=x^3-x^2-2x+1$ as the minimal polynomial for $k$.  Again, using Pari, we compute
$$
\mu(\arithm)=\frac{8 \pi \cdot 49^{3/2}\zeta_k(2)}{(4 \pi^2)^3} = \frac{\pi}{21}.
$$
Thus,
\begin{equation}\label{no3}
 4 \pi = \mu(\Gamma) =M \mu(\arithm) =M\frac{ \pi \prod_{\Prime | \Delta(A)}( N(\Prime)-1) }{21}
\end{equation}
If we take $Ram_f(A) = \emptyset$, then M = 84.  So $\mu(\arithm)=  \frac{\pi}{21}$.  Since $Ram_f(A) = \emptyset$, and $k=k_7$, we have $a_2, a_3, a_7 \neq 0$ by Lemma \ref{split}.  Therefore,
$$
\mu(\arithm) =\frac{ \pi}{21} = 2 \pi \left( 2g-2 +\frac{a_2}{2} + \frac{2 a_3}{3}+  \frac{6 a_7}{3}  \right) 
$$
The only solution to this equation is $a_2=a_3=a_7=1$, and in this case $\arithm$ is a triangle group of signature $(0;2,3,7)$ (see \cite{T1}).  Again the existence of a torsion-free subgroup of $\arithm$ of index 84 is guaranteed by Theorem \ref{torfree}.

If $Ram_f(A) \neq \emptyset$, then $\prod_{\Prime| D(A)} (N(\Prime)-1) \geq 42$.  This is because the primes of smallest norm in $R_k$, $\Prime_2$ and $ \Prime_7$, in $k$ have norms 8 and 7, respectively, and $|Ram_f(A)| \geq 2$.  This implies $M \leq 2$. However, since $\arithm$ will not be torsion-free, Theorem \ref{torfree} implies that $M \geq 2$.  Thus, the only other possible solution to (\ref{no3}) occurs when $M=2$ . 

In this case, $Ram_f (A)= \{ \Prime_2,\Prime_7 \}$ is a solution to (\ref{no3}) when $M=2$.  Since $6 \not| d_k$, we can apply Lemma \ref{sixdivides}.  By Pari, we find that $h=h(k(i))=h(k(\omega))=1$ and that $\Prime_2$ ramifies and $\Prime_7$ is inert in $k(i)|k$; therefore $a_2=2$.  Since the ideal $\Prime_7$ splits in $k(\omega)|k$, we have $a_3=0$.  But
$$
\mu(\arithm) = 42  \cdot \frac{\pi}{21} = 2 \pi \left( 2g-2 + 1 + \frac{6 e_7}{7} \right)
$$
and $g=1, e_7=0$ is the only solution. Thus, $\arithm$ has signature $(1;2,2)$ and, again by Theorem \ref{torfree}, it has a torsion-free subgroup of genus two and of index 2.  Since $k|\rat$ is Galois, there is only one conjugacy class of groups $\arithm$ of this signature.  
\end{proof}  
Johansson \cite{J1} has determined the signatures of all DAFGs arising from quaternion algebras over number fields of degree less than or equal to 2, so we will be concerned with fields of degree $|k:\rat|>2$.  Combining Lemma \ref{list} and our results with the relevant results in \cite{J1},\cite{MR2}, and \cite{T1}, we obtain the following theorem. 

\begin{theorem} \label{list2}
The following is a complete list of all $\arithm$ containing a derived arithmetic Fuchsian group $\Gamma$ of signature $(2;-)$ arising from quaternion algebras over totally real number fields.  The number $c$ denotes the number of conjugacy classes of the group $\arithm$ in each case.
\begin{center}
\begin{tabular}{c||c|c|c|c|c}
%\cline{1-6}
$[k: \rat]$ & {\bf $d_k$} & {\bf $\Delta(A)$}   & $|\arithm:\Gamma|$ & $\arithm$ & $c$ \\\hline\hline
$1$ & $1$& $2\cdot3$ & $6$ & $(0;2,2,3,3)$&$1$\\
$1$ &$ 1$& $2\cdot5$ & $3$ & $(0;3,3,3,3)$&$1$ \\
$1$ &$ 1$& $2\cdot7$ & $2$ & $(1;2,2)$  &$1$\\
$1 $& $1$& $2\cdot13$ & $1$ & $(2;-) $ &$1$ \\\hline
$2$ & $5$ & $\Prime_{2}$ & $20$ & $(0;2,5,5)$& $ 1$  \\
$2$ & $5$ & $\Prime_{5}$ & $15$ & $(0;3,3,5)$ & $1$\\
$2$ & $5$ & $\Prime_{11}$ & $6$ & $(0;2,2,3,3)$ &$ 1$ \\
$2$ &$ 5$ & $\Prime_{11}'$ & $6$ & $(0;2,2,3,3)$ &$1$\\
$2$ & $5$ & $\Prime_{31}$ & $2$ & $(1;2,2)$ &$1$ \\
$2$ & $5$  & $\Prime_{31}'$ & $2$ & $(1;2,2)$ &$1$ \\
$2$ & $5$  & $\Prime_{61}$ & $1 $&$ (2;-)$ & $1$\\
$2$ & $5$  & $\Prime_{61}'$ & $1$ & $(2;-)$ & $1$\\
$2$ & $8$ & $\Prime_{2}$ & $24 $&$ (0;3,3,4)$& $1$ \\
$2$ & $8$ & $\Prime_{5}$ & $1$ & $(2;-)$ & $1$\\

\end{tabular}

\begin{tabular}{c||c|c|c|c|c}
%\cline{1-6}
$[k: \rat]$ & {\bf $d_k$} & {\bf $\Delta(A)$}   & $|\arithm:\Gamma|$ & $\arithm$ & $c$ \\\hline\hline
$2$ & $12$ & $\Prime_{2}$ & $12$ & $(0;3,3,6)$& $1 $ \\
$2$ & $12$ & $\Prime_{3}$ & $6$ & $(0;2,2,2,6)$& $1$  \\
$2$ & $12$ & $\Prime_{13}$ & $1$ & $(2;-)$&$1$   \\
$2$ & $12$ & $\Prime_{13}'$ & $1$ & $(2;-)$ &$1$ \\
$2$ & $13$ & $\Prime_{13}$ & $1 $& $(2;-)$ &$1$\\
$2$ & $17$ & $\Prime_2$ & 6 & $(0;2,2,3,3)$&$1$  \\
$2$ & $17$ & $\Prime_2'$ & 6 & $(0;2,2,3,3)$ &$1$\\
$2$ & $24$ & $\Prime_3$ & $2$ & $(0;2,2,2,2,2,2)$&$1$  \\
$2$ & $28$ & $\Prime_2$ & $3$ & $(0;3,3,3,3)$&$1$  \\\hline
$3$ & $49$ & $\Prime_2 \Prime_7$&  $2$& $(1;2,2)$ & $1$ \\ 
$3$ & $49$ & $\emptyset$ & $84$ & $(0;2,3,7)$ & $1$  \\
$3$ & $81$ & $\emptyset$ & $36$ & $(0;2,3,9)$  &$1$ \\ 
$3$ & $148$ & $\Prime_2 \Prime_{13}$& $1$ & $(2;-)$& $3$   \\
$3$ & $148$ & $\Prime_2 \Prime_5$ &$ 3$ & $(0;3,3,3,3)$&$3$ \\
$3$ & $148$ & $\emptyset$ & $12$ & $(0;2,2,2,3)$& $3$ \\
$3$ & $169$ & $\emptyset$ & $12$ & $(0;2,2,2,3)$&$1$ \\
$3 $& $229$ & $\emptyset$ & $6$ &$(0;2,2,3,3)$&$4$ \\
$3 $& $229$ & $\Prime_2$,$\Prime'_2$& $2$&$(1;2,2)$& $3$\\
$3$ & $257$ & $\emptyset$ & $6$ & $(0;2,2,3,3)$ &$4$ \\
$3$ & $316$ & $\Prime_2$,$ \Prime'_2$ & $3$ &$(0;3,3,3,3)$ &$3$  \\
$4$ & $725$ & $\Prime_2$ & $4$ & $(1;2)$ & $2$ \\
$4$ & $725$& $\Prime_{11}$ & $6$ &$(0;2,2,3,3)$&$2$ \\
$4$ & $725$& $\Prime'_{11}$ & $6$ & $(0;2,2,3,3)$&$2$ \\
$4$ & $725$& $\Prime_{31}$ & $2$ & $(1;2,2)$& $2$\\
$4$ & $725$ & $\Prime'_{31}$ & $2$ & $(1;2,2)$ & $2$\\
$4$ & $725$& $\Prime_{61}$ &$1$ & $(2;-)$ &$2$ \\
$4$ & $725$& $\Prime'_{61}$ &$1$ & $(2;-)$ & $2$\\
$4$ & $1125$ & $\Prime_2$ & $2$ & $(1;2,2)$& $1$ \\
$4$ & $1957$ & $\Prime_7$ & $2$ & $(1;2,2)$& $4$\\ 
$4$ & $1957$ & $\Prime_3$ & $6$ &$(0;2,2,3,3)$&$4$  \\ 
$4$ & $2000$ & $\Prime_5$  & $2$ & $(0;3,3,3,3)$&$2$\\
$4$ & $2304$ & $\Prime_3$  &  $1$ & $(2;-)$& $1$\\
$4$ & $2777$ & $\Prime_2$ & $6$ & $(0;2,2,3,3)$&$1$  \\ 
$4$ & $3981$ & $\Prime_3$  & $2$ &  $(0;2,2,2,2,2,2)$& $4$  \\
$4$ & $4352$ & $\Prime_2$ & $6$ & $(0;3,3,3,3)$&$1$ \\ 
$4$ & $4752$ & $\Prime_2$ & $1$ & $(2;-)$& $2$ \\\hline
$5$ & $24217$ & $\emptyset$ & $12$ & $(0;2,2,2,3)$&$5$\\ 
$5$ & $36497$ &  $\emptyset$ & $6$ & $(0;2,2,3,3)$ &$6$ \\
$5$ & $38569$ &  $\emptyset$ & $6$ &$(0;2,2,3,3)$ &$6$\\
\end{tabular}
\end{center}
\end{theorem}

\begin{remark}   Using a theorem of Greenberg \cite{Gr} on maximal Fuchsian groups in conjunction with the results in \cite{ANR}, \cite{MR2}, \cite{T1}, and \cite{T2} gives all the conjugacy classes of the groups $\arithm$ listed above except for those with signatures $(1;2,2)$, $(0;2,2,2,2,2,2)$, and $(2;-)$.  %In many cases, the number $h_\infty=1$ (see \ref{hinf}).  %We include these in the above Proposition and determine the type number for the three special signatures above.
\end{remark}

\section{Maximal Orders and Fundamental Domains}\label{max}

The group $SU(1,1)$ is the group of orientation-preserving isometries of the unit disk $\unit = \{ z \in \comp | |z| \leq 1 \}$.  By embedding a cocompact arithmetic Fuchsian group $\Gamma$ into $SU(1,1)$, one can determine a fundamental domain for its image, $\Gamma'$, using a theorem of Ford.  The elements of $\Gamma'$ which give the side-pairings of the fundamental domain are generators for $\Gamma'$.  This technique is described for the rational numbers and quadratic number fields in \cite{V} and \cite{K} and more generally  in \cite{J2}.  In order to find a fundamental domain for $\Gamma$ using this technique, the maximal order must be written explicitly as an $R$-module, where $R=R_k$ is the ring of integers of the number field $k$.  We first state and prove some results on the existence and form of maximal orders in certain cases in which the Hilbert symbol for a quaternion algebra $A$ is ``nice".  The invariant quaternion algebras of arithmetic Fuchsian groups with small genus will often fall into this class.     %These will be used in the final chapter to obtain generators for the all derived arithmetic groups of signatures (2;-),(1;2,2), and (0;2,2,2,2,2,2) given in Theorem \ref{list2}.

\subsection{Maximal Orders}\label{maxord}

Recall that any quaternion algebra $A$ has an associated \emph{Hilbert symbol}
$$
\left(  \frac{a,b}{k} \right)
$$
where $i^2 = a, j^2 = b,  ij = -ji = k$ for some $a$, $b$ in $k^*$.  The basis $\{1, i,j,ij \}$ is referred to as the standard basis of $A$.  The discriminant $\Delta(A)$ of a quaternion algebra $A$ is defined to be the product of the prime ideals at which $A$ is ramified.   For any $R$-order $\ord$ in $A$, the discriminant $d(\ord)$ is defined to be the $R$-ideal generated by the set $\{\mbox{det}(\mbox{tr}(x_i x_j)),1 \leq i,j \leq 4 \}$, where $x_i \in \ord$.  We will use the following facts about orders (cf. \cite{MRe}):
\begin{enumerate}
\item[(i)]  Any order is contained in a maximal order.
\item[(ii)]  An order $\ord$ is maximal if and only if $d(\ord)=\Delta (A)^2 $.  %In particular, all maximal orders have the same discriminant.
\item[(iii)]  \label{disc} If $\ord$ has free $R$-basis $\{e_1,e_2,e_3,e_4 \}$, then the discriminant of $\ord$, $d(\ord)$ is the principal ideal $\mbox{det}(\mbox{tr}(e_i e_j))R$.
\end{enumerate}

\begin{proposition}\label{denom}Suppose that $k$ has class number 1 and that $ab$ is square-free where $a,b \in R$.  Let $A=\left(  \frac{a,b}{k} \right)$ be a quaternion algebra over a number field $k$.  Suppose, in addition, that $\Delta(A)$ divides $abR$.   Let  $\pi_i R=\Prime_i$ for each $\Prime_i \not\in Ram_f(A)$ and 
$$r=
\left\{
\begin{array}{cl}
  1&  \quad \Delta(A)=abR   \\
\displaystyle \prod_{\begin{array}{c}  \Prime_i | abR\\ \Prime_i \not| \Delta(A) \end{array}} \pi_i & \quad \Delta(A)\neq abR 
\end{array}
\right..
$$
If $\ord=\ordtwo[\beta]$, where $\ordtwo=R[1,i,j,ij]$, is a maximal order of $A$ for some $\beta \in A$, then 
$$
\beta \in \frac{1}{2r}\ordtwo.
$$ 
\end{proposition}
\begin{proof}  By assumption, all ideals of $R$ are principal.  Since $d(\ordtwo)=16 a^2 b^2 R$, the order $\ordtwo$ is not maximal and $\ordtwo \subset \ord$ for some maximal order $\ord$.  In particular, $d(\ord) = \Delta(A)^2$. Now, since $ab$ is square-free, $abR = r_1 \Delta(A)$, up to multiplication by a unit.     Moreover, $R$ is a principal ideal domain, so both $\ordtwo$ and $\ord$ have free bases, say $\{e_i\}_{i=1}^4$ and $\{f_j \}_{j=1}^4$, respectively.  Since $\ordtwo \subset \ord$, we can write each $e_i$ as $\sum_{j=1}^4 a_{ij} f_j$, where $a_{ij} \in R.$  We therefore have:
$$16a^2 b^2 =d(\ordtwo)=\mbox{det}(\mbox{tr}(e_i e_j)) =(\mbox{det}M)^2 \mbox{det}(\mbox{tr}(f_i f_j))= (\mbox{det}M)^2 d(\ord),$$ 
where $M=(a_{ij})$.  Thus, up to multiplication by a unit, $\mbox{det}M = 4r.$  This implies $\ord \subset \frac{1}{4r}\ordtwo$, i.e., $\beta = \frac{1}{4r}(x_0 + x_1 i + x_2 j + x_3 ij)$, where $x_i \in R$, $0 \leq i \leq 3$.   But since $\beta$ is an integer, its trace must be integral: $\mbox{tr}(\beta)=\frac{x_0}{2r} \in R$.   From this it follows that $x_0 \in 2R.$  Similarly, taking the products $i \beta$, $j \beta$ and $ij \beta$ and using the hypothesis that $\ordtwo[\beta]$ is an order, it follows that $x_1, x_2, x_3 \in 2R.$  Therefore, $\beta \in \frac{1}{2r}\ordtwo$, as claimed.
 
\end{proof}

\begin{lemma}  \label{solmod}  Let $A={\underline{a,b} \choose k}$ be a quaternion algebra over a number field $k$ and ring of integers $R$ with $a,b \in R$ satisfying
  
\begin{enumerate}  
\item[(i)]  $(a,b)=1$,
\end{enumerate}
and either of the following conditions:
\begin{enumerate}  
\item[(ii)] $\exists \, \tilde{a}, \tilde{b} \in R$ such that $\tilde{a}^2 \equiv a \bmod 4R$ and $\tilde{b}^2 \equiv b \bmod 4R$,
\item[(ii)$'$] $b =-1$ and $\exists \, \tilde{a} \in R$ such that $\tilde{a}^2 \equiv a \bmod 4R$.
\end{enumerate}
Then there exists a nonzero solution  $(x,y)\in R \times R$ to the equation  $x^2-ay^2 \equiv b \bmod 4R$.
\end{lemma}

\begin{proof} We need to show that, under the hypotheses, there exist $x, y \in R$ such that $x^2-ay^2 \equiv b \bmod 4R$ or, equivalently, such that $x^2-ay^2-b \in 4R$.  Suppose conditions (i) and (ii) hold.  Then $x^2-ay^2-b \equiv 0 \bmod 4R$ is equivalent to $x^2- \tilde{a}^2 y^2\equiv \tilde{b}^2  \bmod 4R$.  Now, the equation
\begin{equation}\label{glosol}
x- \tilde{a} y \equiv \tilde{b}  \bmod 2R
\end{equation}
will have a solution $(x,y)\in R\times R$ provided 
\begin{equation}\label{losol}
x- \tilde{a}_{\Prime} y \equiv \tilde{b}_{\Prime}  \bmod R_{\Prime}
\end{equation}
has a nonzero solution $( \tilde{x}_{\Prime}, \tilde{y}_{\Prime})$ for every prime $\Prime$ dividing 2 in $R$ (by the Chinese Remainder Theorem).  
Since $(a,b)=1$ implies $( \tilde{a}_{\Prime}, \tilde{b}_{\Prime})=1$, equation (\ref{losol}) clearly has a nonzero solution $(\tilde{x}_{\Prime}, \tilde{y}_{\Prime})$ for each prime $\Prime$ dividing 2.  Therefore, (\ref{glosol}) has a nontrivial solution $(x,y) \in R$.  Since $x- \tilde{a} y \equiv \tilde{b}  \bmod 2R$ if and only if $x+ \tilde{a} y \equiv \tilde{b}  \bmod 2R$, $(x,y)$ satisfies
$$
(x- \tilde{a} y)(x+ \tilde{a} y) \equiv x^2 -\tilde{a}^2y^2  \equiv \tilde{b}^2 \bmod 4R. 
$$
If conditions (i) and (ii)$'$ hold, then the equation $x^2-ay^2-b \equiv x^2 - \tilde{a}^2y^2 - b\equiv 0 \bmod 4R$ is equivalent to
$$
-x^2+\tilde{a}y^2 \equiv 1 \bmod 4R.
$$ 
Again, by the Chinese Remainder Theorem, there exists $(x,y) \in R \times R$ such that $-x+ \tilde{a} y \equiv 1  \bmod 2R$, and the proof now follows as above.

\end{proof}

\begin{lemma} \label{intord} Let $A=\left(\frac{a,b}{k} \right)$ be a quaternion algebra with $a,b \in R$ such that $abR=\Delta(A)$.  Let $\ordtwo=R[1,i,j,ij]$ so that $\ordtwo$ is an order in $A$.  If $\beta \in \frac{1}{2}\ordtwo \setminus \ordtwo$ has the form $\frac{1}{2}(x_0+x_1 i + x_2 j)$ and is integral, then $\ord=R[1,i,\beta,i \beta]$ is a ring of integers.  If, in addition, $x_2 \in R^{*}$, then $\ord \supset \ordtwo$ is a maximal order of $A$.   
\end{lemma}

\begin{proof} Let $e_0=1, e_1=i, e_2=\beta$, and $e_3=i \beta $.  Now $\ord=R[1,i,\beta,i \beta]$ is an order if and only if the following conditions are satisfied:
\begin{enumerate}
\item[(i)] $ e_k e_l$ is integral for $0 \leq k,l \leq 3$
\item[(ii)] $e_k + e_l$ is integral for $0 \leq k,l \leq 3.$
\end{enumerate}
The simple structure of this order makes many of these conditions redundant.  The conditions in (i) and (ii) are conditions that the norms and traces of these elements belong to $R$.  Moreover, (i) and (ii) also establish that $\ord$ is closed under multiplication.  The norms and traces of these elements are listed in Tables 1 and 2 below.  Note that although the elements $e_k e_l$ and $e_l e_k$, $k \neq l$, may not be equal, their traces and norms are equal (hence, both tables are symmetric).  We have also omitted the obvious cases, e.g., 1 and $i$.

Since $a, b, x_k \in R$, for $0 \leq k \leq 2$, all of the conditions on integrality reduce to the following conditions:
\begin{enumerate}
\item[(i)] $x_0^2-ax_1^2-bx_2^2 \in 4R$
\item[(ii)] $(x_0^2-ax_1^2-bx_2^2)^2 \in 16R$
\item[(iii)] $x_0^2+ax_1^2+bx_2^2 \in 2R.$
\end{enumerate}

Condition (i) implies the all the others.  We will show (i) implies (iii).  The condition $x_0^2-ax_1^2-bx_2^2 \in 4R \Rightarrow x_0^2-ax_1^2-bx_2^2 \in 2R$ since $4R \subset 2R$.  But $x_0^2-ax_1^2-bx_2^2 \equiv x_0^2-ax_1^2-bx_2^2 \bmod 2R$, so $x_0^2-ax_1^2-bx_2^2 \in 2R \Leftrightarrow x_0^2+ax_1^2+bx_2^2 \in 2R$.  However, condition (i) is equivalent to the integrality of $\beta$.  This shows that integrality of $\beta$ implies the integrality of all elements of $\ord=R[1,i,\beta,i \beta]$.  Thus, if $\beta$ is integral, then $\ord$ is a ring of integers.

\begin{table}
\begin{center}
\begin{tabular}{c||c|c|c|c} \label{sums}
%\cline{1-5}
$\times$ & 1 & $i$ & $\beta$ & $i \beta$ \\\hline\hline
1&*&*& n=$\frac{(x_0^2-ax_1^2-bx_2^2)}{4}$& n=$-\frac{a(x_0^2-ax_1^2-bx_2^2)}{4}$\\
& &&$\mbox{tr}=x_0$&$\mbox{tr}=a x_1 $\\\hline
$i$ &*&*&*&n=$\frac{a^2(x_0^2-ax_1^2-bx_2^2)}{4}$\\
&&&&$\mbox{tr}=ax_0$\\\hline
$\beta$&*&*&n=$\frac{(x_0^2-ax_1^2-bx_2^2)^2}{16}$&n=$\frac{a(x_0^2-ax_1^2-bx_2^2)^2}{16}$\\
&&&$\mbox{tr}=\frac{(x_0^2+ax_1^2+bx_2^2)}{2}$&$\mbox{tr}=ax_0 x_1 $\\\hline
 $i\beta$ &*&*&*&n=$-a(x_0^2-ax_1^2-bx_2^2)$\\
 &&&&$\mbox{tr}=2ax_1$\\
\end{tabular}
\end{center}
\caption{Norms and traces of sums of the $R$-basis of $\ord$ in Lemma \ref{intord}}
\end{table}

\begin{table} \label{prods}
\begin{center}
\begin{tabular}{c||c|c|c|c} 
%\cline{1-5}
+ & 1 & $i$ & $\beta$ & $i \beta$ \\\hline\hline
1&*&*& n=$\frac{(4+4x_0-x_0^2-ax_1^2-bx_2^2)}{4}$& n=$\frac{(4+4ax_1-a(x_0^2-ax_1^2-bx_2^2))}{4}$\\
& &&$\mbox{tr}=2+x_0$&$\mbox{tr}=2+a x_1 $\\\hline
$i$ &*&*&n=$\frac{(-4a+4ax_1+x_0^2-ax_1^2-bx_2^2)}{4}$&n=$-\frac{a(4+4x_0+x_0^2-ax_1^2-bx_2^2)}{4}$\\
&&&$\mbox{tr}=x_0$&$\mbox{tr}=ax_1$\\\hline
$\beta$&*&*&n=$(x_0^2-ax_1^2-bx_2^2)$&n=$\frac{(a-1)(x_0^2-ax_1^2-bx_2^2)}{4}$\\
&&&$\mbox{tr}=2x_0$&$\mbox{tr}=x_0+a x_1 $\\\hline
 $i\beta$ &*&*&*&n=$-a(x_0^2-ax_1^2-bx_2^2)$\\
 &&&&$\mbox{tr}=2ax_1$\\
\end{tabular}
\end{center}
\caption{Norms and traces of products of the $R$-basis of $\ord$ in Lemma \ref{intord}}
\end{table}

We will now assume that $x_2 \in R^{*}$.  In order to show that $\ord$ is an order, we must show that $R[1,i,\beta,i \beta]$ is a complete R-lattice with 1.  It is clear that $\ord$ is an $R$-lattice.  Since $1, i \in \ord$, it remains to show that $j \in \ord$.  Since $\beta=\frac{1}{2}(x_0+x_1 i+x_2 j) \in \ord$, we have $j=x_2^{-1}(2\beta-x_1-x_2i) \in \ord$ and, hence, $I$ is complete. Therefore, $  R[1,i,\beta,i \beta]$ is an order.  Moreover, the discriminant of the order $R[1,i,\beta,i\beta]$ is $d(\ord)=a^2b^2x_2^4R=a^2b^2R=\Delta A^2$, since $x_2 \in \units$.   Hence, $\ord=R[1,i,\beta,i \beta]$ is maximal. 
\end{proof}

\begin{proposition} \label{nice} Let $A=\left(\frac{a,b}{k} \right)$ be a quaternion algebra over a number field $k$ with finite ramification set $Ram_f(A)$ and denote the standard order of $A$ by $\ordtwo=R[1,i,j,ij]$.  Suppose that $a, b \in R_k$ satisfy the hypotheses of Lemma \ref{solmod} and, in addition, that $\Delta(A) =abR$.  Then there exists $\beta \in \frac{1}{2}\ordtwo$ such that $\ordtwo \subset \ord$ and $\ord=R[1,i,\beta,i\beta]$ is a maximal order of $A$. 
\end{proposition}
\begin{proof}  As noted previously, the order $\ordtwo=R[1,i,j,k]$ has discriminant $d(\ordtwo)=16a^2b^2R$ and is therefore not maximal.  Since any order is contained in a maximal order, $\ordtwo \subset \ord$ where $\ord$ is a maximal order.  In particular, $d(\ord)=a^2b^2R =\Delta(A)^2$. The ideal $I=\frac{1}{2}\ordtwo \supset \ordtwo$ is not an order since its elements, e.g., $1/2$, are not all integral.  But the discriminant of $I$ is equal to $a^2b^2R$.  Therefore, $\ordtwo \subset \ord$.  This implies there exists some $\beta \in \frac{1}{2}\ordtwo$ such that $\beta \in \ord$.  

Since $a,b$ satisfy the hypothesis of Lemma \ref{solmod}, there exist integers $x_0, x_1 \in R$ such that $x_0^2-ax_1^2-b \in 4R$.  Therefore, if we take $\beta =\frac{1}{2}(x_0 + x_1 i + j) \in \frac{1}{2}\ordtwo$, then $\beta$ is integral.  Furthermore, by Lemma \ref{intord}, $I = R[1,i,\beta,i \beta]$ is a maximal order.  
\end{proof}
If the Hilbert symbol of $A$ does not satisfy the conditions of the previous proposition, we can still use Proposition \ref{denom} as a starting point, but the process of finding a free-basis for $\ord$ becomes more ad hoc.  One uses an intermediate order $\ordthree$ where $\ordtwo \subset \ordthree \subset \frac{1}{2}\ordtwo$ with $d(\ordthree)=a^2b^2R$ and searches for integral elements in the ideal $\frac{1}{r}\ordthree$ where $r \in R$ is as stated in the proof of Lemma \ref{denom}.   By testing these integral elements $\beta$ as part of a free $R$-basis of the orders $R[1,i,\beta,i \beta]$  and $R[1,j, \beta, j \beta]$ and computing the discriminants of these orders, one can determine a maximal order in the algebra.

\subsection{Fundamental Domains and Generators}

Let $A$ be the invariant quaternion algebra corresponding to arithmetic Fuchsian group $\arithm$.  For any maximal order $\ord$ of $A$, fix an embedding $\rho$ of ${\ord}^1$ in $PSL_2(\real)$ and denote the image by $\arithm$.  Choose $\rho$ so that $i \in \hyp$ is not the fixed point of any nontrivial element in $\arithm$.  The Mobius transformation
$$
\varphi = \left( \begin{array}{cc} i & 1 \\
1 & i
\end{array}
\right)
$$    
maps $\hyp$ to the unit disk $\unit$.  Furthermore, the action of $SL_2(\real)$ on $\hyp$ is conjugate to the action of $SU(1,1)$ on $\unit$ since

$$  
SU(1,1) = \varphi SL_2(\real) \varphi^{-1}.
$$
This defines an embedding of $\arithm$ in $SU(1,1)$.  

For any $g \in SU(1,1)$ or $SL_2(\real)$, 
$$ 
g = \left( \begin{array}{cc} a & b \\
c & d
\end{array}
\right)
$$
with $c \neq 0$, the isometric circle $C_g$ of $g$ is defined to be the set of points on which $g$ acts as a Euclidean isometry.  %By taking the derivative of $g$, one can verify that the circle $\{z \in \comp | |cz+d|=1\}=\{z \in \comp | |z+\frac{d}{c}|=\frac{1}{|c|}\}$ is the isometric circle of $g$.    

The following theorem of Ford (cf. \cite{K}, Ch. 3) characterizes a fundamental domain of $\Gamma \subset SU(1,1)$ in terms of the isometric circles of its elements:

\begin{theorem} \label{Ford} Let $\Gamma$ be a discrete subgroup of $SU(1,1)$ such that the origin is not a fixed point of any nontrivial element of $\Gamma$.  Let $C_g$ be the isometric circle of $g$.  If $C_g^o$ is the set of all points outside $C_g$, then
$$
\fund = \unit \cap \bigcap_{g \in \Gamma} C_g^o
$$
is a fundamental domain of $\Gamma$.
\end{theorem}

Clearly, $\varphi^{-1}(\fund)$ is a fundamental domain for $\varphi^{-1}(\Gamma)\varphi$.   Let $r_g$ be the radius of the isometric circle $C_g$ where $g \in \Gamma$ for a discrete subgroup $\Gamma$ of $SU(1,1)$.  Since 
$$
SU(1,1) =\left\{ \left( \begin{array}{cc} a & c \\
\bar{c} & \bar{a}
\end{array}
\right) \Big| a, c \in \comp,\, a\bar{a}-c\bar{c}=1 \right\},
$$
the radius $r_g = \frac{1}{|c|}$.  From the discreteness of  $\Gamma$ and the additional relation $a\bar{a}-c\bar{c}=1$,  it follows that 
$$
\Gamma_{\epsilon}=\{g \in \Gamma | r_g > \epsilon \} 
$$
is finite for every $\epsilon$, $0 <\epsilon< 1$.  If, for some $ \epsilon>0$, 
$$
\fund_{\epsilon} = \unit \cap \bigcap_{g \in \Gamma_{\epsilon}} C_g^o,\quad\unit{\epsilon}=\{z \in \comp | |z| < 1- \epsilon\}, \quad \fund_{\epsilon}\subset \unit_{\epsilon},
$$
then $\fund_{\epsilon}$ will be a fundamental domain for $\Gamma$.  This will be the case for some sufficiently small $\epsilon >0$, since $\Gamma$ has finite coarea and no parabolic elements.  

Using this consequence of Ford's Theorem one can systematically obtain the generators for arithmetic Fuchsian groups if one can find free $R$-basis for the maximal order.   We use this technique in conjunction with our results on maximal orders to obtain generators for some of the unit groups $\arithm$ listed in Theorem \ref{list2}.   Our main interest will be in the cases $3 \leq |k:\rat| \leq 4$, since examples of this type are lacking in the literature.  Although this technique is described in generality in \cite{J2}, we will include a description here for completeness.

Let $r_g$ denote the radius of the isometric circle $\varphi g \varphi^{-1}$, where $g \in \Gamma \subset PSL_2(\real)$.  If
$$ 
g= \left( \begin{array}{cc} a & b \\
c & d
\end{array}
\right),
$$
then
$$
\varphi g \varphi^{-1} = \frac{1}{2}\left( \begin{array}{cc} (a+d)+i(b-c) & (b+c)+i(a-d) \\
(b+c)-i(a-d)  & (a+d)-i(b-c) 
\end{array}
\right).
$$
Therefore,
\begin{equation} \label{radius}
r_g = \frac{2}{|(b+c)-i(a-d)|}=\frac{2}{\sqrt{(a-d)^2+(b+c)^2}}=\frac{2}{\sqrt{a^2+b^2+c^2+d^2-2}}.
\end{equation}
Hence, the restriction $r_g > \epsilon$ gives an upper bound on the entries of $g$.   Furthermore, using the fact that the norm is positive definite in all other $n-1$ embeddings of $\sigma_i:k \hookrightarrow \rat$, one obtains upper bounds on the absolute values of $\sigma_i(a),\sigma_i(b),\sigma_i(c),\sigma_i(d)$ for $2\leq i \leq n$.   

If we write $\ord$ as a $\ints$-module, then we use the bounds on the $\sigma_i$ to get bounds on the integral coefficients of the elements of $\Gamma$.  

Let $|k: \rat| =n$ and suppose that $k$ has integral power basis $\{1, \al, \ldots, \al^{n-1}\}$.  By properties of quaternion algebras over the real numbers, we may assume that $A$= $\underline{a,b} \choose k $, where $a>0$ and $b<0$.  Fix an embedding $\rho: A \hookrightarrow M_2(k (\sqrt{a}))$.  Then the standard order $R[1,i,j,ij]$ is the set of elements
$$
\left(\begin{array}{cc} x +y\sqrt{a}& b_1(u+v\sqrt{a}) \\
b_2(u-v\sqrt{a}) & x-y\sqrt{a}
\end{array}
\right)
$$
where $b_1, b_2 \in R$ satisfy $b_1 b_2 = b$.  Now, by Proposition \ref{denom}, a maximal order $\ord$ of $A$ is contained in $\frac{1}{r}R[1,i,j,ij]$, for some $r \in R \setminus \{ 0\}$.  Therefore, $\ord$ will be a subset of the set of elements of the form
\begin{equation} \label{intvector}
g=\frac{1}{r}\left(\begin{array}{cc} (\sum x_i \al^i)+ (\sum y_i \al^i)\sqrt{a}& b_1( (\sum u_i \al^i)+ (\sum v_i \al^i)\sqrt{a}) \\
b_2( (\sum u_i \al^i)- (\sum v_i \al^i)\sqrt{a}) &  (\sum x_i \al^i)- (\sum y_i \al^i)\sqrt{a}
\end{array}
\right),
\end{equation}
where $r \in R \setminus \{ 0\}$ and the integers $x_i, y_i, u_i, v_i \in \ints$, $0 \leq i \leq n-1$.  The integrality of the elements in $\ord$ translate to certain congruence relations on the $x_i, y_i, u_i, v_i \in \ints$, $0 \leq i \leq n-1$.  The norm of $g$ is $\mbox{n}(g)= \frac{1}{r^2}(x^2-ay^2-bu^2+abv^2)$.    
Since norm is invariant under each embedding $\sigma_i$ of the number field, $\mbox{n}(g)=1$ implies $\mbox{n}(\sigma_i({g}))=1$, for $1 \leq i \leq n.$  Therefore, for each $g \in \arithm$, we have
$$
\sigma(r)^2=\sigma(x)^2-\sigma(a)\sigma(y)^2-\sigma(b)\sigma(u)^2+\sigma(a)\sigma(b)\sigma(v)^2.
$$
Since $A$ is a quaternion algebra ramified at all but one finite place, we may assume that $\sigma_1(a) > 0 $, $\sigma_1(b) <0 $, $\sigma_i(a) <0$, and $\sigma_i(b) <0$ for $2 \leq i \leq n$.  Therefore, for each i, $2 \leq i \leq n$,
\begin{equation}\label{bounds1}
\begin{array}{lll} 
|\sigma_i(x)| &\leq  &  \sigma(r),  \\  
|\sigma_i(y)| & \leq &\sqrt{\frac{\sigma(r)^2-\sigma_i(x)^2}{-\sigma_i(a)}},\\
|\sigma_i(v)| & \leq &\sqrt{\frac{\sigma(r)^2-\sigma_i(x)^2+\sigma_i(a)\sigma_i(y)^2}{\sigma_i(a)\sigma_i(b)}},\\ 
|\sigma_i(u)| & \leq & \sqrt{\frac{\sigma(r)^2-\sigma_i(x)^2+\sigma_i(a)\sigma_i(y)^2-\sigma_i(a)\sigma_i(b)\sigma_i(v)^2}{-\sigma_i(b)}}. \\
\end{array}
\end{equation}

Substituting into (\ref{radius}), we get that $r_g=\frac{2}{\sqrt{q-2}}$, where
$$
q=\frac{1}{r^2}\left( 2x^2 + 2ay^2 +\frac{4a b^2}{b_1^2+b_2^2}v^2 +(b_1^2+b_2^2)(u+v\frac{b_1^2-b_2^2}{b_1^2+b_2^2} \sqrt{a})^2 \right).
$$
The condition that $r_g > \epsilon$ is equivalent to $q < M_{\epsilon}:= 2 + \frac{4}{\epsilon^2}$, and this condition implies the following set of bounds for $\sigma_1=\mbox{Id}$:
\begin{equation}\label{bounds2}
\begin{array}{lll} 
|x| & < &  r \sqrt{\frac{M_{\epsilon}}{2}},  \\  
|y| & < &   \sqrt{\frac{1}{2a}(r^2 M_{\epsilon} - 2x^2)},  \\  
|v| & < &   \sqrt{\frac{b_1^2+b_2^2}{4ab^2}(r^2 M_{\epsilon} - 2x^2-2ay^2)},  \\  
|u| & < &   \sqrt{\frac{1}{b_1^2+b_2^2}(r^2 M_{\epsilon} - 2x^2-2ay^2-\frac{4ab^2}{b_1^2+b_2^2}v^2)}+\frac{b_1^2-b_2^2}{b_1^2+b_2^2}\sqrt{a}|v|.  \\  
\end{array}
\end{equation}

By taking various linear combinations of these inequalities, we obtain bounds on the integers  $x_i, y_i, u_i, v_i \in \ints$, $0 \leq i \leq n-1.$  

\section{Examples}

In this section, we use our previous results to find generators for a few examples of the DAFGs $\arithm$ in Theorem \ref{list2} with signature $(1;2,2)$, $(0;2,2,2,2,2,2)$, or $(2;-)$ using programs written in Mathematica (see \cite{M}, Appendix B).  (A complete list of generators for all the groups $\arithm$ with one of these signatures can be found in \cite{M}, Ch. 5).  Using a standard presentation of the group $\arithm$ and Magma, we also explicitly determine generators for each subgroup $\Gamma$ of signature $(2;-)$. %Example programs for number fields of degrees 2 and 3 are given in the appendix.  

Here we give examples in which the Hilbert symbol of $A$ satisfies the hypotheses of Proposition \ref{nice} and examples in which it does not.  All elements of $\arithm$ will be given as a vector of integers using an integral power basis of $R$ with a specified denominator $r$ (cf. (\ref{intvector})).  Also, we will abuse notation and use the same vector to describe the corresponding matrix in $SL_2(\real)$.  

\begin{example}\label{nicecubic}  Let $k=\rat( \cos(\frac{\pi}{7}))=k(\al)$ be the totally real cubic field of degree 3, where $\al$ is a root of the polynomial $f(x)=x^3-x^2-2x+1$.  The group $\arithm$ with invariant quaternion algebra $A$ defined over $k$ with $Ram_f(A) =\{\Prime_2,\Prime_7 \}$   %A Hilbert symbol for $A$ is ${ \underline{2(2\al-3),-1} \choose k}$ 
is generated by the elements 
\begin{center}
 \begin{tabular}{c||c|c|c|c|c|c|c|c|c|c|c|c} 
% \cline{1-13}
 & $x_1$& $x_2$ & $x_3$ & $y_1$& $y_2$ & $y_3$&$u_1$& $u_2$ & $u_3$& $v_1$& $v_2$ & 
 $v_3$\\\hline\hline
$A_1$ &0&2&2&-1&2&2&0&2&2&1&-2&-2 \\
$B_1$ &-3&2&3&1&-2&-2&1&-2&-3&-2&2&3 \\
$X_1$ &0&0&0&0&0&0&2&0&0&0&0&0 \\
\end{tabular}\\
\end{center}
where $r=2$ and $<A_1,B_1, X_1|([ A_1,B_1]X_1)^2,X_1^2>$. Here the vectors $A_1$, $B_1$, and $X_1$ are as described at (\ref{intvector}) with $a = 2(2 \al -3)$ and $b=-1$.
\end{example}
%\begin{eqnarray*}
%a_1&=& \frac{1}{2}\left(-4+5\alpha+7\alpha^2+(-1+\alpha^2)i +(\alpha+\alpha^2)j+(3-5\alpha-6\alpha^2)ij \right)\\
%b_1&=&\frac{1}{2}\left( \alpha+\alpha^2+(3-5\alpha-6\alpha^2)i + (4-5\alpha-7\alpha^2)j +(-1+\alpha^2)ij \right)\\
%c_1&=& j
%\end{eqnarray*}
%where $i^{2}=2(2\alpha-3), j^{2}=-1$.

\begin{figure}[h] \label{cubic49}
\centerline{
\epsfig{file=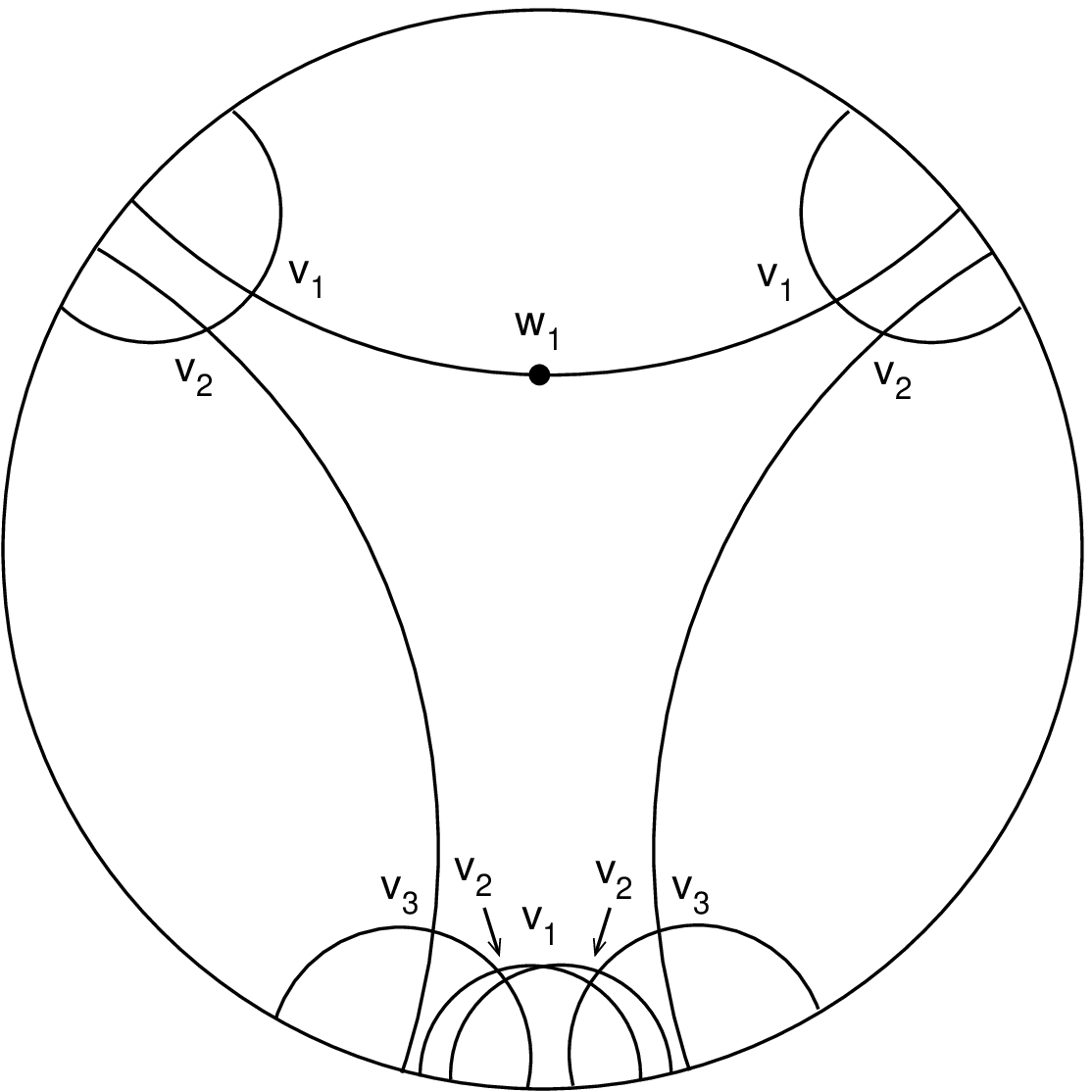, scale=0.6}
}
\caption{Fundamental region for $\arithm$ of Example \ref{nicecubic}}
\end{figure} 

%\begin{figure}[h] \label{cubic49}
%\centerline{
%\scalebox{0.5}{\pdfimage{cubic49.pdf}}
%}
%\caption{Fundamental region for $\arithm$ in Example \ref{nicecubic}}
%\end{figure} 

\begin{proof}
The cubic field $k_7$ has minimal polynomial $f(x)=x^3-x^2-2x+1$ and discriminant 49.  By Proposition \ref{list2}, there is only one conjugacy class of groups $\arithm$ of signature $(1;2,2)$ defined over $k$.  If we denote the three roots of $f(x)$ by $\alpha_1,\alpha_2$ and $\alpha_3$, where $\alpha_1 < 0 <\alpha_2 <\alpha_3$, then the algebra
$$
\left( \frac{2(2\alpha-3),-1}{k} \right).
$$
has the correct ramification set.  This can easily be checked using standard results in algebra (cf. \cite{MRe}, Ch. 2).   %By Pari, there is one unique prime lying over 2 and one unique prime lying over 7 with $\Prime_7 = (7,\alpha+2)R=(2\alpha-3)R$ and $\Prime_2=2R$.   
%We will first show that the quaternion algebra $A$ is given by the Hilbert symbol Since $2\alpha_1-3,2\alpha_2-3 <0$ and $2\alpha_3-3>0$, $A$ is ramified at the two real places corresponding to $\alpha_1$ and $\alpha_2$.  %Furthermore, $-1 \equiv 6 \bmod \Prime_7$ and $6$ is not a square $\bmod 7$, so $\Prime_7 \in Ram_f(A)$ by Theorem \ref{finiteram}.  By Theorems \ref{finiteram} and \ref{card}, $A$ cannot be ramified at any nondyadic prime other than $\Prime_7$; therefore, $Ram_f(A)=\{\Prime_2, \Prime_7 \}$.  This shows that $A$ has the correct ramification set.

In this case, one can check that $A$ does not satisfy the hypotheses of Proposition \ref{nice}.  However, since $\Prime_2=2R$ and $\Prime_7 = (2\alpha-3)R$, we have $abR= 2(2\alpha-3)R=\Prime_2 \Prime_7 = \Delta(A)$.   Therefore, a maximal order $\ord$ of $A$ will nonetheless be of the form $R[1,i,\beta, i\beta]$ where $\beta \in \frac{1}{2}\ord' \backslash \ord'$, where $\ord'=R[1,i,j,ij]$.  The element $\beta = \frac{1}{2}(1+i+j)$ is integral since $\mbox{tr}(\beta)=1 \in R$ and $\mbox{n}(\beta)=-\al+2 \in R$ and $d(\ord')=\Prime_2^2 \Prime_7^2 = \Delta(A)^2$ so that $\ord'$ is maximal.  We can write $\frac{1}{2}R[1,i,j,ij]$ as the $\ints$-module 

$$\left\{ \frac{1}{2}\left(\sum_{k=1}^3 x_k \al^k + i \sum_{k=1}^3  y_k \al^k   +j \sum_{k=1}^3  u_k \al^k  +ij \sum_{k=1}^3 v_k \al^k \right) \Big| x_i, y_i, u_i, v_i  \in \ints \right\}.$$  

Similarly, we write $\ord=  (m_i   + n_i i + o_i \beta + p_i i \beta)$, where $m_i, n_i, o_i, p_i \in \ints[1, \al, \al^2, \al^3]$.  Clearly, $\ord \subset \frac{1}{2}R[1,i,j,ij]$.  Setting the two $\ints$-modules equal yields a linear system of equations.  Since the $m_i,n_i,o_i,p_i$ are integers, solving the system for these variables gives congruence conditions on the $x_i,y_i, u_i, v_i$.  These are necessary and sufficient conditions for an element $g \in  \frac{1}{2}R[1,i,j,ij]$ to be an element of $\ord$.  In this particular case, $\ord$ is the set of elements of $\frac{1}{2}R[1,i,j,ij]$ satisfying the following congruence relations:
$$
\begin{array}{rl}
x_0 + u_0 \equiv 0 \bmod 2 \quad&\quad y_0+u_0+v_0 \equiv 0 \bmod 2\\
x_1+u_1 \equiv 0 \bmod 2 \quad&\quad y_1+u_1+v_1 \equiv 0 \bmod 2\\
x_2+u_2 \equiv 0 \bmod 2 \quad&\quad y_2+u_2+v_2 \equiv 0 \bmod 2  \\
\end{array}.
$$

%In the setting of (\ref{intvector}), 
We implement the inequalities (\ref{bounds1}) and (\ref{bounds2}) with the values $r=2, b_1=2$, and $b_2=-1/2$ and, in this case, $\epsilon =.15$ is sufficient to obtain the Ford domain for $\arithm$.

The Ford domain for the group $\arithm$ is shown in Figure 1.   Since $v_1$ and $v_2$ are distinct fixed points of elements of order two, $\arithm$ has signature $(1;2,2)$.  The elements listed in the Table 3 are the generators for $\arithm$ corresponding to the side pairings of the Ford domain.

 \begin{table}\label{tabnicecubic}
 \begin{center}
 \begin{tabular}{c||c|c|c|c|c|c|c|c|c|c|c|c} 
% \cline{1-13}
 & $x_1$& $x_2$ & $x_3$ & $y_1$& $y_2$ & $y_3$&$u_1$& $u_2$ & $u_3$& $v_1$& $v_2$ & $v_3$\\\hline\hline
$ h_1$&0&2&2&1&-2&-2&0&2&2&1&-2&-2\\ 
$ h_2$&0&2&2&1&-2&-2&0&-2&-2&-1&2&2\\ 
$ h_3$&-1&0&1&0&0&0&-1&0&1&1&0&-1\\ 
$ h_4$&-1&0&1&0&0&0&-1&4&5&3&-4&-5\\ 
$ g_1$&0&0&0&0&0&0&2&0&0&0&0&0\\ 
\end{tabular}
\end{center}
\caption{Generators for $\arithm$ in Example \ref{nicecubic}}
\end{table}

The elements $A_1= h_2^{-1}$, $ h=h_1 h_3^{-1},  X_1=g_1$ are noncommuting hyperbolic elements that satisfy the relation $([A_1,B_1] X_1)^2 = -I$.  Furthermore, no proper subrelation is trivial.  Therefore, if we denote the group $<A_1, B_1, X_1|([A_1,B_1] X_1)^2,X_1^2>$by $\Gamma'$, then $\Gamma'$ has signature $(1;2,2)$.  Since $h_1= A_1^{-1} X_1$ and $ h_4= A_1^{-1} B_1^{-1}$, we have $<A_1, B_1, X_1> \subset \arithm$.  But since Fuchsian groups are Hopfian, $\arithm$ cannot contain a proper isomorphic subgroup.  Therefore, $\displaystyle \arithm = <A_1, B_1, X_1|([ A_1, B_1]X_1)^2, X_1^2>$ and $ A_1, B_1, X_1$ are generators for $\arithm$.
%We can use Magma to determine generators for all subgroups of index two in a group of this signature.  In this case, $\arithm$ contains the three distinct torsion-free subgroups of index two given in the statement of the proposition.  

In this case, one can easily check that the group $\arithm$  has four distinct subgroups $\Gamma_i$, $1 \leq i \leq 4$, of signature $(2;-)$ of index two.  Using the standard presentation $$\displaystyle <A_1, B_1, X_1, Y_1| [[ A_1, B_1]X_1Y_1, X_1^2, Y_1^2]>$$ 
for $\arithm$, we use Magma to find generators for all the subgroups of $\arithm$ of index 2.  Of these, there are four subgroups that are torsion-free, which we will denote by $\Gamma_i$, $1 \leq i \leq 4$.  The presentations for each subgroup are:
$$
\begin{array}{ll}
\Gamma_1 = < B_1A_1^{-1},  X_1A_1^{-1},Y_1A_1^{-1}, A_1^{-2} > &  \Gamma_3 = < A_1,  X_1B_1^{-1}, Y_1B_1^{-1}, B_1^{-2}>\\
\Gamma_2 = < B_1, X_1A_1^{-1}, Y_1A_1^{-1},  A_1^{-2}> &\Gamma_4 = <A_1, B_1, X_1A_1X_1, X_1B_1X_1>
\end{array}.
$$

For each $\Gamma_i$, we determine the trivial relation in the group.  After putting each group in the standard presentation $\displaystyle < a_i, b_i,c_i,d_i |[a_i,b_i][c_i, d_i] >, 1 \leq i \leq 4,$ 
we obtain the corresponding list of generators (see Table 4).

 \begin{table} \label{genustwocub}
 \begin{center}
 \begin{tabular}{c||c|c|c|c|c|c|c|c|c|c|c|c} 
 %\cline{1-13}
 & $x_1$& $x_2$ & $x_3$ & $y_1$& $y_2$ & $y_3$&$u_1$& $u_2$ & $u_3$& $v_1$& $v_2$ & $v_3$\\\hline\hline
$ a_1$&-14&20&24&-12&18&22&-12&20&24&12&-18&-22 \\
$ b_1$&-15&20&25&12&-18&-22&15&-20&-25&-13&18&23 \\
$ c_1$&0&2&2&-1&2&2&0&2&2&-1&2&2 \\
$ d_1$&30&-42&-52&-30&42&53&-32&46&58&-26&38&47 \\\hline

%& &  & & & & & &  & & &  & \\\hline
$ a_2$&-3&2&3&1&-2&-2&1&-2&-3&-2&2&3 \\
$ b_2$&-14&20&24&-12&18&22&-12&20&24&12&-18&-22 \\
$ c_2$&30&-42&-52&-30&42&53&-32&46&58&-26&38&47 \\
$ d_2$&0&2&2&-1&2&2&0&2&2&-1&2&2 \\\hline

$ a_3$&0&2&2&-1&2&2&0&2&2&1&-2&-2 \\
$ b_3$&-14&21&25&13&-18&-23&18&-25&-31&-15&23&28 \\
$ c_3$&1&-2&-3&-2&2&3&3&-2&-3&-1&2&2\\ 
$ d_3$&-37&52&65&-33&47&59&-3&2&3&2&-3&-4\\\hline
$ a_4$&0&2&2&-1&2&2&0&2&2&1&-2&-2 \\
$ b_4$&-3&2&3&1&-2&-2&1&-2&-3&-2&2&3 \\
$ c_4$&0&-2&-2&-1&2&2&0&-2&-2&1&-2&-2 \\
$ d_4$&3&-2&-3&1&-2&-2&-1&2&3&-2&2&3 \\
\end{tabular}
\end{center}
\caption{Generators for $\Gamma_i, 1 \leq i \leq 4$ in Example \ref{nicecubic}}
\end{table}
\end{proof}

\begin{example}\label{nicequartic}  Let $k=\rat(\al)$ be the totally real quartic field of discriminant 3981 where $\al$ is a root of the polynomial $f(x)=x^4-x^3-4x^2+2x+1$.  The group $\arithm$ corresponding to the quaternion algebra $A$ with $Ram_f(A)=\{ \Prime_3 \}$ and that is unramified at the infinite place corresponding to the root $\alpha_2$, where $-1 < \al_2 < 0$, defined over $k$ is generated by 
\begin{center}
 \begin{tabular}{c||c|c|c|c|c|c|c|c|c|c|c|c|c} 
 %\cline{1-14}
 & $x$  & $y_1$ & $y_2$ & $y_3$ & $y_4$ & $u_1$ & $u_2$ & $u_3$ & $u_4$ & $v_1$ & $v_2$ & $v_3$ & $v_4$ \\\hline\hline
$ X_1$ &0& -2 & 4 & 1 & -1& 2& -4 & -1&1&-3&4 &1 &-1  \\
$ X_2$  &0&0 & 0 &0 & 2&0&0&0&0&0&0&0&0  \\
$ X_3$  &0& -2 & 4 & 1 & -1& -2& 4 & 1&-1&3&-4 &-1 &1  \\  
$ X_4$ &0& -3 & 4 & 1 & -1& -17& 18 & 7&-5&36&-39 &-15 &11    \\
$ X_5$  &0& 0 & 0 & 0 & 0 & 20 & -22 & -8& 6& -42 & 50 & 18 & -14  \\
\end{tabular}
\end{center}
where $r= 2$ and $ \displaystyle < X_1, X_2, X_3, X_4, X_5 | X_1^2, X_2^2,  X_3^2, X_4^2, X_5^2, ( X_1 X_2 X_3 X_4 X_5)^2>.$  Here $a = -\al(\al+1)$ and $b = -1$.
\end{example}

\begin{proof}

Let $\al_1 < -1 < \al_2 < 0 < \al_3 < 1 < \al_4$ denote the four roots of $f(x)$.  %In this case, there are four conjugacy classes of groups $\arithm$ defined over $k$.  
The algebra $A=~{ \underline{-\al(\al+1),-1} \choose k}$ is unramified at the place $\sigma_2$ since $-\al_i(\al_i + 1) < 0$ for $i=1,3,4$ and $-\al_2(\al_2 + 2) >0$.  There is a unique prime of norm 3 in $R_k$: $\Prime_3 = (3,\al+1)R=(\al+1)R$.  Furthermore, one can easily verify that  %$f(x) \equiv (x+1)(x^2+1) \bmod 3$, so there are two primes dividing 3.  Let $\Prime_3 = (3,\al+1)R=(\al+1)R$.  Since $-1 \equiv 2 \bmod \Prime_3 \equiv 2 \bmod 3$ and 2 is not a square $\bmod 3$, $A$ is ramified at $\Prime_3$ by Theorem \ref{finiteram}.  Since no other prime divides $ab$ and there is a unique prime above 2, by Theorems \ref{finiteram} and \ref{card}, 
$Ram_f(A)=\{ \Prime_3 \}$.  This is a ``nice" Hilbert symbol, so Proposition \ref{nice} applies, and we find that $\ord = R[1,i,\beta, i \beta]$ where $\beta = \frac{1}{2}(1+\al + \al^2 i + j)$ is a maximal order.  The congruence relations in this case are:
$$
\begin{array}{rr}
x_0 + u_0 +u_3 + v_0\equiv 0 \bmod 2 \quad&\quad y_0+u_2+u_3 +v_0+v_3 \equiv 0 \bmod 2\\
x_1+ u_0 +v_1\equiv 0 \bmod 2 \quad&\quad y_1+ u_3 +v_0+v_1 \equiv 0 \bmod 2\\
x_2+u_1 +u_2+v_2\equiv 0 \bmod 2 \quad &\quad y_2+ u_0+v_1 +v_2\equiv 0 \bmod 2 \\
x_3+ u_2+v_3\equiv 0 \bmod 2 \quad& \quad y_3+ u_1+u_2 + u_3 +v_2\equiv 0 \bmod 2.  \\
\end{array}
$$

\begin{figure}[h] \label{quartic3981}
\centerline{
\epsfig{file=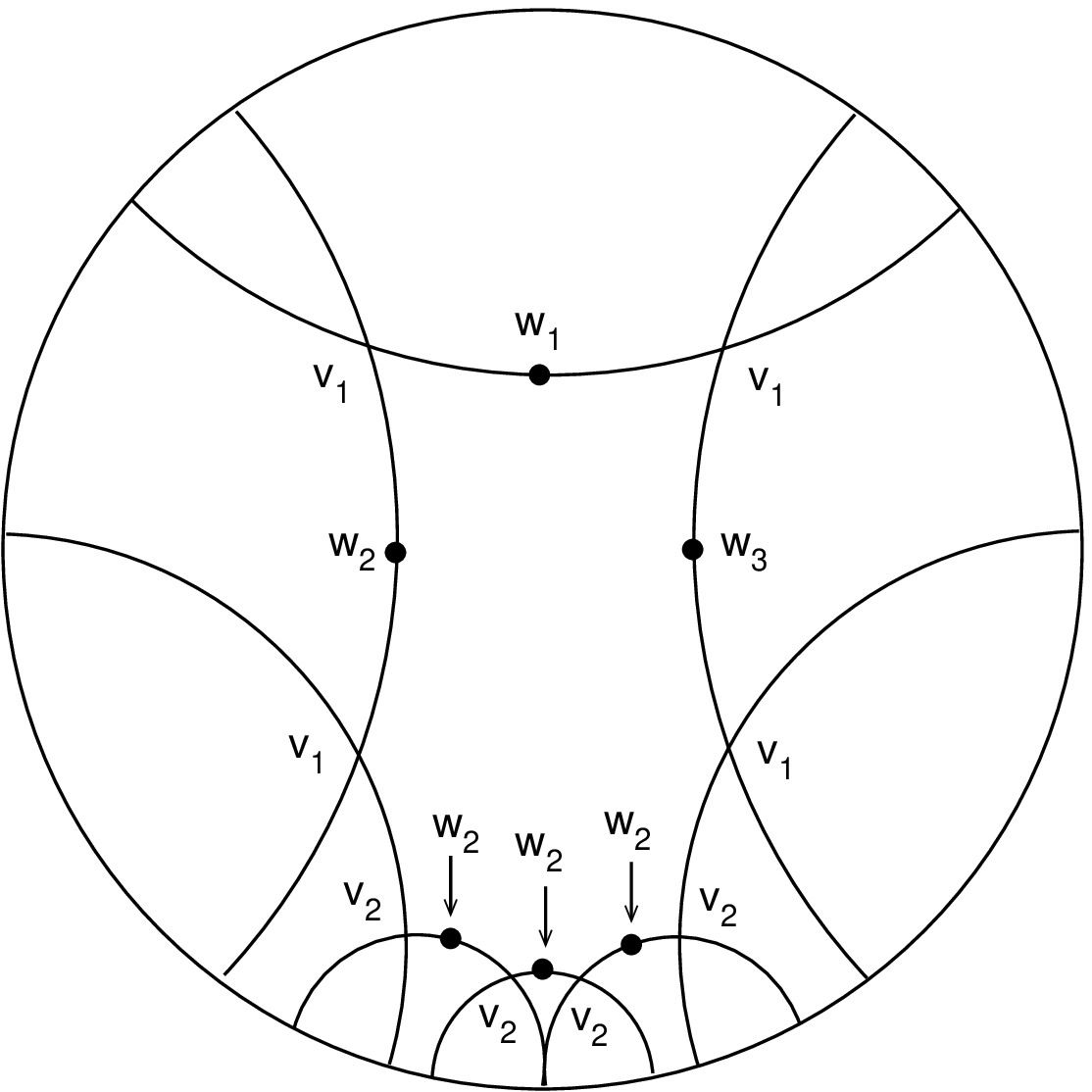, scale=0.6}
}
\caption{Fundamental region for $\arithm$ of Example \ref{nicequartic}}
\end{figure}

In this case, we implement the inequalities (\ref{bounds1}) and (\ref{bounds2}) using $r=2, b_1=2, b_2=-1/2,$ and $\epsilon = .15$. 
We obtain the fundamental region shown in Figure 2 and the corresponding generators are listed in Table 5.  The points $w_i$ are the fixed points of the $g_i$, $1 \leq i \leq 6$,  which all have order two; this verifies that $\arithm$ has signature (0;2,2,2,2,2,2).

 \begin{table}\label{genustwoquart}
 \begin{center}
 \begin{tabular}{c||c|c|c|c|c|c|c|c|c|c|c|c|c|c|c|c} 
 %\cline{1-17}
$\!\!$ &$\!\! x_1\!\!$ $\! $  & $\!\! x_2\!\!$ $\!\!$ & $\!\! x_3 \!\!$  & $\!\! x_4 \!\!$ & $\!y_1\!$& $\!y_2\!$ & $\!\!y_3\!\!$&$\!\!y_4\!\!$&  $\! \!u_1\!\!$ & $\!\! u_2\!\!$ &   $\!\!u_3\!\!$ & $\!\!u_4\!\!$ &  $\! \!v_1\!\!$  & $\!\! v_2 \!\!$ & $\!\!v_3\!\!$ & $\!\!v_4\!\!$\\\hline\hline

$\! g_1\!$  &$\!\!$  0 $\!\!$  & $\!\!$  0  $\!\!$ & $\!\!$  0  $\!\!$ &  $\!\! $  0 $\!\!$  & $\!0 \!$&  $\!$ 0 $\!$&  $\!$ 0$\!$ &$\!\!$ 0 $\!\!$& $\!\!$ -2 $\!\!$ & $\!\!$ 0 $\!\!$ & $\!$ 0$\!$ &$\!\!$ 0$\!\!$&  $\!\!$ 0$\!\!$ &  $\!\!$ 0$\!\!$ &  $\!\!$ 0$\!\!$ &  $\!\!$ 0$\!\!$ \\
$\! g_2\!$  & $\!\!$  0 $\!\!$  & $\!\!$ 0  $\!\!$ & $\!\!$  0  $\!\!$ &  $\!\!$  0 $\!\!$ & $\!-2 \!$&  $\!$ 4 $\!$&  $\!$ 1$\!$ &$\!\!$ -1$\!\!$ & $\!\!$ -2$\!\!$& $\!\!$ 4$\!\!$&$\!$ 1$\!$&$\!$ -1$\!$&  $\!$ $\!\!3 \!\!$&  $\!$ $\!\!$-4$\!\!$ &  $\!\!$ -1$\!\!$ &  $\!\!$ 1 $\!\!$\\
$\! g_3\!$  & $\!\!$  0 $\!\!$  & $\!\!$  0  $\!\!$ & $\!\!$  0  $\!\!$ &  $\!\!$  0 $\!\!$  & $\!-2 \!$&  $\!$ 4 $\!$&  $\!$ 1$\!$ &$\!\!$ -1$\!\!$ & $\!\!$ 2$\!\!$& $\!\!$ -4$\!\!$&$\!$ -1$\!$& $\!$1$\!$ &  $\!\!$ -3 $\!\!$&  $\!\!$ 4 $\!\!$&  $\!\!$ 1$\!\!$ &  $\!\!$ -1$\!\!$ \\
$\!g_4\!$  & $\!\!$  0 $\!\!$  & $\!\!$  0  $\!\!$ & $\!\!$  0  $\!\!$ &  $\!\!$  0 $\!\!$  & $\!-3\!$ &  $\!$ 4$\!$ &  $\!$ 1$\!$ &$\!\!$ -1$\!\!$ & $\!\!$ -17$\!\!$&  $\!\!$ 18$\!\!$&$\!$ 7$\!$& $\! $-5$\!$ &  $\!\!$ 36 $\!\!$&  $\!\!$ -39 $\!\!$&  $\!\!$ -15$\!\!$ &  $\!\!$ 11$\!\!$\\
$\!g_5\!$  & $\!\!$   0 $\!\!$  & $\!\!$  0  $\!\!$ & $\!\!$  0  $\!\!$ & $\!\!$   0 $\!\!$  & $\!-3 \!$&  $\!$ 4 $\!$&  $\!$ 1$\!$ &$\!\!$ -1$\!\!$ & $\!\!$ 17$\!\!$&  $\!\!$ -18$\!\!$ &$ \!$ -7$\!$ & $\!5 \!$ &  $\!\!$ -36$\!\!$ &  $\!\!$ 39 $\!\!$&  $\!\!$ 15 $\!\!$&  $\!\!$ -11$\!\!$\\
$\!g_6\!$  & $\!\!$  0 $\!\!$  & $\!\!$  0  $\!\!$ & $\!\!$  0  $\!\!$ &  $\!\!$  0 $\!\!$  & $\!0 \!$&  $\!$ 0$\!$ &  $\!$ 0$\!$ &$\!\!$ 0$\!\!$ & $\!\!$ 20$\!\!$& $\!\!$ -22$\!\!$&$\!$ -8$\!$&$\! 6 \!$&  $\!\!$ -42 $\!\!$&  $\!\!$ 50 $\!\!$&  $\!\!$ 18 $\!\!$&  $\!\!$ -14$\!\!$ \\
$\! h_1\!$  & $\!\!$  2 $\!\!$  & $\!\!$  -3  $\!\!$ &  $\!\!$ -1  $\!\!$ &  $\!\!$  1 $\!\!$  &$\!0\!$ &  $\!$ 0$\!$ &  $\!$ 0$\!$ &$\!\!$ 0$\!\!$ &  $\!\!$ 5 $\!\!$& $\!\!$ -3$\!\!$&$\!$-2$\!$& $\!$1$\!$ &  $\!\!$ -10$\!\!$&  $\!\!$ 11$\!\!$& $\!\!$ 4$\!\!$ &  $\!\!$ -3$\!\!$ \\
\end{tabular}
\end{center}
\caption{Generators for $\arithm$ of Example \ref{nicequartic}}
\end{table}

After putting the group in the required presentation, we obtain the list as stated above.  A group of  signature $(0;2,2,2,2,2,2)$  has a unique subgroup $\Gamma$ of signature $(2;-)$ by \cite{Gr}.  Using Magma, we find that, if $\arithm$ is presented in the form
$$
 <X_2 X_1, X_3 X_1, X_4 X_1, X_5 X_1 |  X_1^{-1} X_2 X_3^{-1}X_4 X_1 X_2^{-1}X_2 X_4^{-1}>,
$$
then the subgroup $\Gamma$ is generated by $< X_2X_1, X_3X_1, X_4X_1, X_5X_1>$.  Putting these generators in the standard form $< a_1,b_1,c_1,d_1|[a_1,b_1][c_1,d_1]>$ yields the set of generators for $\Gamma$ listed in Table 6.

\begin{table}
\begin{center}
 \begin{tabular}{c||c|c|c|c|c|c|c|c|c|c|c|c|c|c|c|c} 
 %\cline{1-17}
$\!\!$ & $\!x_1\!$  &$x_2$  &$x_3$  &$x_4$  & $y_1$ & $y_2$ & $y_3$ & $y_4$& $u_1$ & $u_2$ & $u_3$ & $u_4$ & $v_1$ & $v_2$ & $v_3$ & $v_4$ \\\hline\hline
$\!\!a_1\!\!$& $\!9\!$&-11&-4&3&67&-75&-28&31&-35&39&15&-11&39&-43&-16&12\\
$\!\!b_1\!\!$&$\!2\!$&-4&-1&1&-10&11&4&-3&7&-7&-3&2&-11&14&5&-4  \\
$\!\!c_1\!\!$&$\!$-12$\!$&14&5&-4&26&-28&-11&8 &-14&18&6&-5&31&-36&-13&10\\
$\!\!d_1\!\!$&$\!$-2$\!$&4&1&-1&31&-36&-13&10&-14&18&6&-5&3&-4&-1&1 \\
\end{tabular}
\end{center}
\caption{Generators for $\Gamma$ in Example \ref{nicequartic}}
\end{table}\end{proof}

%\begin{tabular}{|c|c|c|c|c|c|c|c|c|} 
% \cline{1-9}
% & $u_1$ & $u_2$ & $u_3$ & $u_4$ & $v_1$ & $v_2$ & $v_3$ & $v_4$ \\\hline
%$a_1$&-35&39&15&-11&39&-43&-16&12 \\\hline
%$b_1$&7&-7&-3&2&-11&14&5&-4 \\\hline
%$a_2$&-14&18&6&-5&31&-36&-13&10 \\\hline
%$b_2$&-14&18&6&-5&3&-4&-1&1 \\\hline
%\end{tabular}

Our last example is defined over a quartic number field in which the Hilbert symbol of A does not satisfy the hypotheses of Proposition \ref{nice}.  We will show that a ``nice" symbol does not exist for $A$. We remark that the technique for finding a maximal order is similar to that of the previous examples, except in this case we cannot take $r$ to be a rational integer.

\begin{example} \label{badquartic}   Let $k = \rat(\al)$ be the totally real quartic field of discriminant 4752, where $\al$ is a root of the polynomial $f(x)=x^4-2x^3-3x^2+4x+1$. The group $\arithm$ of signature $(2;-)$ defined over $k$ with quaternion algebra $A$ satisfying $Ram_f(A)=\{ \Prime_2 \}$ and unramified at the infinite place corresponding to the root $\alpha_1$, where $-2 < \al_1 < -1$  %Then a Hilbert symbol for $A$ is $ {\underline{1+\alpha,-(-1+\alpha)(-1+\alpha+\alpha^{2})} \choose \rat(\alpha)}$ and, if $\ord \subset A$ is a maximal order, then $\arithm$ 
has the following list of generators:
\begin{center}
 \begin{tabular}{c||c|c|c|c|c|c|c|c|c|c|c|c|c|c|c|c} 
% \cline{1-17}
$\!$ $\!$& $\! x_1\!$ & $\!x_2\!$ &  $\!x_3\!$ & $\!x_4\!$ & $\!y_1\!$ & $\!y_2\!$ & $\!y_3\!$ & $\!y_4\!$ & $\!u_1\!$ & $\!u_2\!$ & $\!u_3\!$ & $\!u_4\!$ & $\!v_1\!$ & $\!v_2\!$ & $\!v_3\!$ & $\!v_4\!$ \\\hline\hline
$\!  h_1\!$ &$\!1\!$ & 0 & -1 & 0 & -1 & 0 & 0 & 0& 0& 0 & 0&0&0&4 &1 &-1  \\
$\! h_2\!$ &$\!0 \!$& -2 & -1 & 1 & -1 & 3 & 1 & -1&0&0&0&0&-1&-2&1&0  \\
$\! h_3\!$ & $\!0\!$& -4 & 0 & 1 & 1 & -1 & 0 & 0 & 1 & 3 & 1 & -1& 1 & 2 & -1 & 0  \\
$\! h_4\!$& $\!0\!$& -4 & 0 & 1 & 1 & -1 & 0 & 0 & -1 & -3 & -1 & 1& 1 & 2 & -1 & 0  \\
$\! h_5\!$ & $\!1\!$& 0 & -1 & 0 & 0 & 3 & -1 & 0 & 1 & 1 & 2 & -1& 0 & -1 & -2 & 1  \\
$\! h_6\!$ & $\!1\!$& 0 & -1 & 0 & 0 & 3 & -1 & 0 & -1 & -1 & -2 & 1& 0 & -1 & -2 & 1  \\
$\! h_7\!$ & $\!1\!$& 2 & -2 & 0 & 2 & 2 & -1 & 0 & 0 & 0 & 0 & 0& -1 & -1 & 3 & -1  \\
$\! h_8\!$ & $\!0\!$& -1 & -3 & 1 & -1 & 3 & -1 & 0 & 0 & 0 & 0 & 0& 0 & -1 & 3 & -1  \\
$\! h_9\!$ & $\!1 \!$& 1 & -2 & 1 & 0 & 3 & -1 & 0 & 0 & 0 & 0 & 0& 0 & -5 & 2 & 0  \\
\end{tabular}
\end{center}
%\begin{eqnarray*}
%a_{1}&=&  \frac{1}{2(1+\alpha)}\left(4\alpha-\alpha^{3}+(-1+\alpha)i +(-1-3\alpha-\alpha^{2} +\alpha^{2})j \right. \\
%&&\quad \left.+(-1-2\alpha+\alpha^{2})ij \right),\\
%b_{1}&=& \frac{1}{1+\alpha}\left(2\alpha-3\alpha^{2}+\alpha^{3}+(1+3\alpha-4\alpha^{2}+\alpha^{3})i \right.\\
%&& \quad \left. +(2+5\alpha-5\alpha^{2}+\alpha^{3})j+(3+9\alpha-9\alpha^{2}+2\alpha^{3})ij  \right), \\
%a_{2}&=&\frac{1}{2(1+\alpha)}\left(6+15\alpha-28\alpha^{2}+8\alpha^{3}+(36+120\alpha-186\alpha^{2}+53\alpha^{3})i \right.\\
%&&\quad \left.+(117+388\alpha-606\alpha^{2}+173\alpha^{3})j \right.\\
%&&\quad \left.+(159+528\alpha-831\alpha^{2}+238\alpha^{3})ij  \right),  \\
%b_{2}&=&\frac{1}{2(\alpha+1)}\left(-1-2\alpha+2\alpha^{2}+(2+12\alpha-15\alpha^{2}+4\alpha^{3})i \right. \\
%&&\quad \left.+(9+28\alpha-50\alpha^{2}+15\alpha^{3})j+(13+40\alpha-68\alpha^{2}+20\alpha^{3})ij \right),
%\end{eqnarray*}
%where $i^{2}=1+\alpha$ and $j^{2}=(1-\alpha)(-1+\alpha+\alpha^{2})$.
Here $a =1+\alpha$ and $b=(1-\alpha)(-1+\alpha+\alpha^{2})$.
\end{example}
\begin{proof}
Let $\alpha_1 < -1 < \alpha_2 < 0 < 1<\alpha_3 < 2 < \alpha_4$ denote the four real roots of $f(x)$.  %Since $f(x) \equiv (x^2 +x +1)^2$, there is a unique prime $\Prime_2$ of norm 4 dividing 2; furthermore, $\Prime_2 = (2,\alpha^2+\alpha+1)=(\alpha^2+\alpha-1)=(\pi)$.  
A fundamental system for $R^{*}$ is $<\alpha,\alpha-1,\alpha^2-2>$.  The signs of the generators of $R^{*}$ and of the uniformizer $\pi$ for $\Prime_2$ under the different embeddings corresponding to the $\alpha_i$ are shown below:
$$
\begin{array}{c|cccc}
& \alpha & -1+\alpha & -2+\alpha^{2} & \pi \\\hline 
\alpha_1 \sim -1.4955& - & - & + & -\\
\alpha_2 \sim -0.21968& - & - & - & -\\
\alpha_3\sim 1.2196 & + & + & - & +\\
\alpha_4 \sim 2.4955& + & + & + & +\\
\end{array} 
$$
From the table of embeddings, it is evident that there does not exist a Hilbert symbol for $A$ such that the only primes dividing $abR$ are in $Ram_f(A)$.  In this case, the algebra
$$
\left(\frac{1+\alpha,(1-\alpha)(-1+\alpha+\alpha^{2})}{\rat(\alpha)} \right)
$$
is unramified at the place $\sigma_1$. The element $1+\alpha$ is  a uniformizer for $\Prime_3$ in R (since $f(x) \equiv (x+1)^4 \bmod 3$, $\Prime_3$ is the unique prime of norm 3 lying over 3).  But $A$ is unramified at $\Prime_3$ since $(1-\alpha)(-1+\alpha+\alpha^2 )\equiv 1 \bmod \Prime_3$ which is a square $\bmod \; 3$.  Thus, $A$ corresponds to one of the two conjugacy classes of $\arithm$ defined over $k$.

In this case, we use an intermediate order to find a maximal order $\ord$.  The order $\ordtwo = R[1,i,\gamma, i\gamma]$, where $\gamma=\frac{1}{2}((1+\alpha)+(1+\alpha^2)i+j)$ has discriminant $d(\ordtwo)= (1+\alpha)^2(-1+\alpha)^2(-1+\alpha+\alpha^2)^2R = \Prime_3^2 \Prime_2^2R \neq \Delta(A)^2$, and therefore is not maximal.  Thus, $\ordtwo \subset \ord$, for some maximal order $\ord.$  Again, by arguments analogous to those in the proof of Prop. \ref{denom}, one can argue that there exists an element $\beta \in \ord \cap \frac{1}{\alpha+1}\ordtwo \backslash \ordtwo$.  An element in $\ordtwo$ written as an integral vector satisfies the following congruence relations:
\begin{equation}\label{intermedord}
\begin{array}{rl}
x_0 +u_0+u_3+v_0 +v_1 +v_2 +v_3& \equiv 0 \bmod 2\\
x_1 +u_0 +v_0+v_1+v_2+v_3&\equiv 0 \bmod 2\\
x_2 +u_1-u_2+u_3+v_0 +v_3  &\equiv 0 \bmod 2\\
x_3 +u_2 +u_3 +v_0+v_1&\equiv 0 \bmod 2 \\
y_0 -u_0+u_2+v_0+v_3 &\equiv 0 \bmod 2 \\
y_1+u_1 +u_3+v_0+v_1&\equiv 0 \bmod 2 \\
y_2 +u_0+v_1+v_2+v_3&\equiv 0  \bmod 2 \\
y_3 +u_1+v_2+v_3&\equiv 0  \bmod 2. \\
\end{array}
\end{equation}
Consider an element of the form $\beta=\frac{1}{2(\alpha+1)}(x+ y i + u j +  v ij) \in I=\frac{1}{\alpha+1}\ordtwo \backslash \ordtwo$.  Now, $\beta$ is integral if and only if 
\begin{eqnarray*}
\mbox{tr}(\beta)&=&\frac{x}{1+\alpha} \in R \\
\mbox{det}(\beta)&=&\frac{1}{4(1+\alpha)^2} \left( x^2+(1+\alpha) y^2+(-1+\alpha)(-1+\alpha+\alpha^{2})u^2 \right.\\
&& \quad \left.  + (1+\alpha)(-1+\alpha)(-1+\alpha+\alpha^{2})v^2 \right)  \\
&=&  \frac{d(x,y,u,v)}{4(1+\alpha)^2} = \frac{d}{4(1+\alpha)^2}. \\
\end{eqnarray*}

If we write each of $x, y,u,v$ in an integral power basis of $R$, these conditions are equivalent to the existence of solutions $(r_0,r_1,r_2,r_3),(s_0,s_1,s_2,s_3) \in \ints^4$ to the equations
\begin{equation}\label{intsol}
\begin{array}{rll}
x  &=& (1+\alpha)(r_0+r_1 \alpha+r_2 \alpha^2 + r_3 \alpha^3 )\\
& =& (r_0-r_3)+(r_0+r_1-4r_3)\alpha+ (r_1+r_2+3r_3)\alpha^2+(r_2+3r_3)\alpha^3\\
d &=& 4(1+\alpha)^2(s_0+s_1 \alpha+s_2 \alpha^2 + s_3 \alpha^3 )\\
&=& 4(s_0 -s_2 - 4s_3)+4(2s_0+s_1-4s_2-17s_3)\alpha\\
&& \quad + 4(s_0+2s_1+4s_2+8s_3)\alpha^2+4(s_1+4s_2+12s_3)\alpha^3 \\
\end{array}
\end{equation}
where $d=d(x_0,\ldots x_4, y_0,\ldots y_4,u_0,\ldots u_4,v_0,\ldots v_4) \in R$.  The expressions for $x$ and $d$ are simplified by the relation $\alpha^4=2\alpha^3+3\alpha^2-4\alpha-1$.  The existence of solutions to these equations yields another set of congruences on the $x_i,y_i,u_i,v_i$.  Using these in addition to the congruence relations (\ref{intermedord}), we find that 
$$
\beta = \frac{1}{2(\alpha+1)}\left((1+\alpha+\alpha^2+\alpha^3)+(1+\alpha+2\alpha^2)i+ij \right) 
$$
is integral; $\mbox{tr}(\beta)=1+\alpha^2 \in R$ and $\mbox{det}(\beta)=2\alpha^3+1$.  Furthermore, 
$$
\begin{array}{rll}
i \beta &=& \frac{1}{2(\alpha+1)} \left((1+\alpha+2\alpha^2)(\alpha+1)+ (1+\alpha+\alpha^2+\alpha^3)i+(\alpha+1)j \right) \\
 &= &\frac{1}{2}\left((1+\alpha+2\alpha^2)+ (1+\alpha^2)i+j \right) \\
\end{array};
$$
this implies $j = 2i \beta - (1+\alpha+2\alpha^2) - (1+\alpha^2)i \in I$.  The integrality of the elements of $I$ are verified using (\ref{intsol}) and the traces and products of the $R$-basis are listed in Tables 7 and 8.  %The solutions  $(r_0,r_1,r_2,r_3),(s_0,s_1,s_2,s_3)$ to (\ref{intsol}) are also given in Table 4.6.

\begin{table} \label{sumsbad}
\begin{center}
\begin{tabular}{c||c|c|c|c} 
%\cline{1-5}
$\times$ & 1 & $i$ & $\beta$ & $i \beta$ \\\hline\hline
1&*&*& n=$3+\alpha^2+2\alpha^3$& n=$-1-8\alpha+4\alpha^2+6\alpha^3$\\
& &&$\mbox{tr}=3+\alpha^2$&$\mbox{tr}=1-\alpha-2\alpha^2$\\\hline
$i$ &*&*&*&n=$1-5\alpha+7\alpha^2+7\alpha^3$\\
&&&&$\mbox{tr}=-1-\alpha-2\alpha^2$ \\\hline
$\beta$&*&*&n=$4(1+2\alpha^3)$& n=$-7\alpha+6\alpha^2+8\alpha^3$\\
&&&$\mbox{tr}=2(1+\alpha^2)$& $\mbox{tr}=-\alpha-\alpha^2$\\\hline
 $i\beta$ &*&*&*&n=$-4(1+7\alpha-6\alpha^2-6\alpha^3)$\\
 &&&&$\mbox{tr}=-2(1+\alpha+2\alpha^2)$
\end{tabular}
\end{center}
\caption{Norms and traces of sums of the $R$-basis of $\ord$ in Example \ref{badquartic}}
\end{table}

\begin{table} \label{prodsbad}
\begin{center}
\begin{tabular}{c||c|c|c|c} 
%\cline{1-5}
+ & 1 & $i$ & $\beta$ & $i \beta$ \\\hline\hline
1&*&*& n=$1+2\alpha^3$& n=$-1-7\alpha+6\alpha^2+6\alpha^3$\\
& &&$\mbox{tr}=1+\alpha^2$&$\mbox{tr}=1+\alpha^2 $\\\hline
$i$ &*&*&n=$-1-7\alpha+6\alpha^2+6\alpha^3$&n=$1-3\alpha+10\alpha^2+8\alpha^3$\\
&&&$\mbox{tr}=-1-\alpha-2\alpha^2$&$\mbox{tr}=-3-3\alpha-2\alpha^2$\\\hline
$\beta$&*&*&n=$(1+2\alpha^3)^2$&n=$-37-170\alpha+70\alpha^2+108\alpha^3$\\
&&&$\mbox{tr}=-2-4\alpha+5\alpha^2-2\alpha^3$&$\mbox{tr}=-3-5\alpha+3\alpha^2-2\alpha^3$\\\hline
 $i\beta$ &*&*&*&n=$(-1-7\alpha+6\alpha^2+6\alpha^3)^2$\\
 &&&&$\mbox{tr}=2(-1+5\alpha^2)$
 \end{tabular}
 \end{center}
\caption{Norms and traces of products of the $R$-basis of $\ord$ in Example \ref{badquartic}}
\end{table}

Since $R[1,i,\beta, i \beta]$ has discriminant $\Prime_2^2 =\Delta(A)^2$, it is a maximal order and, hence, $\ord =R[1,i,\beta, i \beta]$.  Finally, we determine the congruence relations for $\ord$ written as a $\ints$-module:
\begin{equation}
\begin{array}{rl}
4x_0-x_1+x_2-x_3 +u_0+2u_1-2u_2+2u_3-3v_0  +3v_2&\!\!\! \equiv 0 \bmod 6\\
x_0 -x_1+x_2-x_3-2u_0+2u_1-8u_2 -7u_3 -v_1+v_1&\!\!\! \equiv 0 \bmod 6 \\
4y_0 -y_1+y_2-y_3-u_0-2u_1+2u_2+u_3+v_0+2v_1-2v_2+2v_3 &\!\!\! \equiv 0 \bmod 6 \\
y_0-y_1+y_2-y_3+2u_0-2u_1-u_2+u_3 -2v_0+2v_1-2v_2-v_3&\!\!\! \equiv 0  \bmod 6 \\
-u_0+u_1-u_2+u_3 &\!\!\! \equiv 0 \bmod 3 \\
x_2+x_3  +u_1+v_1+v_3&\!\!\! \equiv 0 \bmod 2\\
x_0 +x_1+x_2+u_2+v_0   &\!\!\! \equiv 0 \bmod 2\\
y_2+y_3+u_0 +u_1+v_1&\!\!\! \equiv 0 \bmod 2 \\
y_0 +y_1+y_2+u_1+u_2+u_3+v_2 &\!\!\! \equiv 0  \bmod 2 \\

\end{array}.
\end{equation}
Here we use $b_1= \pm b_2=\pm \sqrt{(\alpha-1)(\alpha^2-\alpha+1)}$ and $\epsilon=0.1$ in (\ref{bounds1}) and (\ref{bounds2}).
The fundamental region for $\arithm$ is shown in Figure 3 and the correponding generators are those listed at the beginning of the example. \end{proof}

\vspace{2.5mm}

\begin{figure}[h] \label{quartic4752}
\centerline{
\epsfig{file=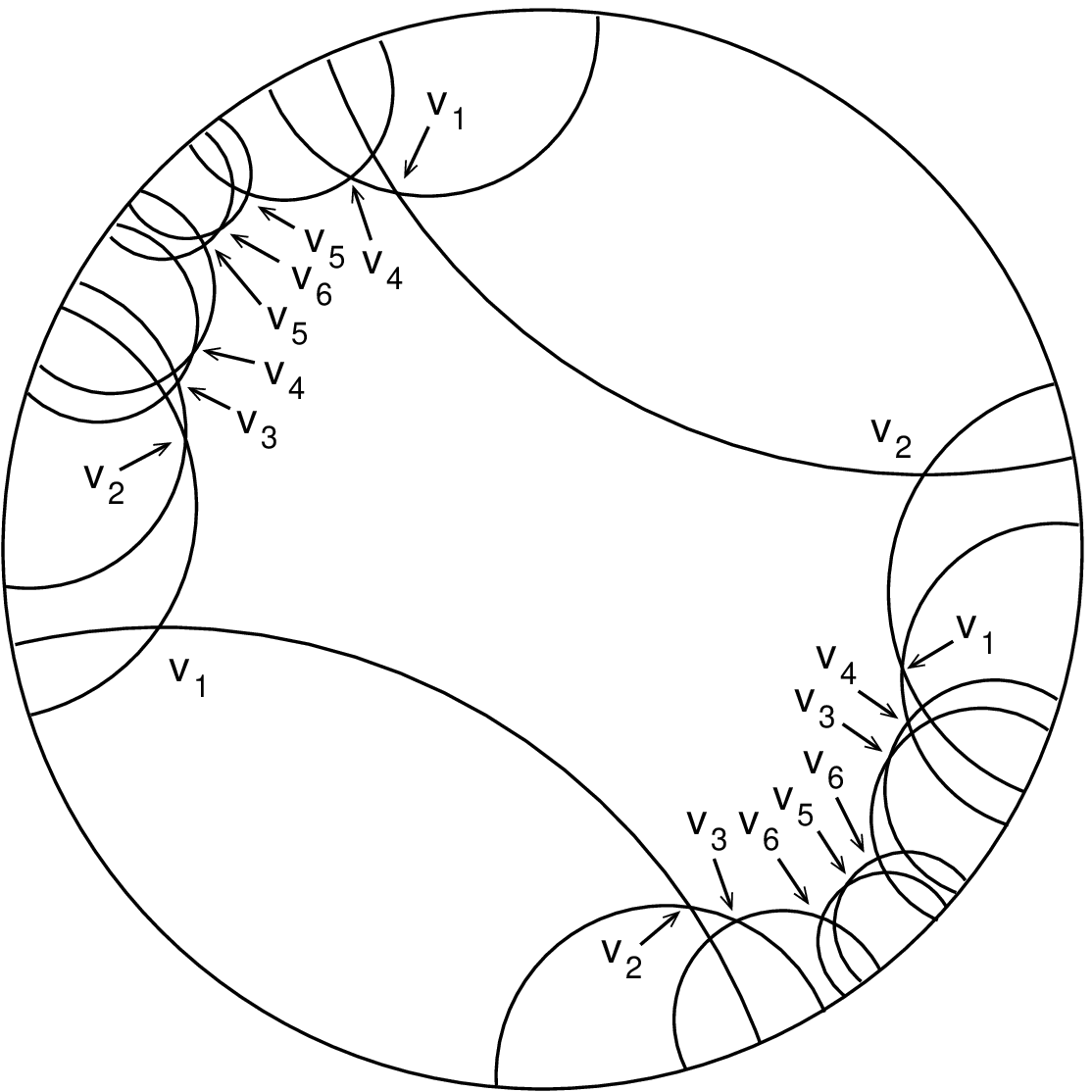, scale=0.6}
}
\caption{Fundamental region for $\arithm$ of Example \ref{badquartic}}
\end{figure} 

%\vbox{
%\begin{center}
%\scalebox{0.5}{\pdfimage{quartic4752.pdf}}\\
%{\sc Figure 3}:  Fundamental region for $\arithm$ of Proposition \ref{badquartic}
%\end{center}
%}
\begin{center}
{\sc Acknowledgements}
\end{center}
This work comprises a large part of the author's Ph.D. thesis and, as such, the author is sincerely grateful and indebted to her advisor, Alan W. Reid, and would like to thank him for suggesting the problem and for the many enlightening conversations concerning it.   %Much of this work would not have been possible without his careful guidance.  
The author would also like to thank Felipe Voloch for pointing out the connection between solutions to the quadratic form in Lemma \ref{solmod} and the Chinese Remainder Theorem and Fernando Rodriguez-Villegas for helpful comments regarding the ramification of quaternion algebras at dyadic primes (cf. \cite{M} Ch. 5).  Lastly, the author would like the referee for the many helpful comments and suggestions contributing to the improvement of this paper. This research was partially supported by an NSF VIGRE fellowship at the University of Texas at Austin.

\nocite{*}      % This command causes all items in the               %
                % bibliographic database to be added to              %
                % the bibliography, even if they are not             %
                % explicitly cited in the text.                      %
                %                                                    %
\bibliographystyle{plain}  % Here the bibliography                   %
\bibliography{diss}

\begin{thebibliography}{10}

\bibitem{ANR}
P.~Ackermann, M.~N{\"a}{\"a}t{\"a}nen, and G.~Rosenberger.
\newblock The arithmetic {F}uchsian groups with signature {$(0;2,2,2,q)$}.
\newblock In {\em Recent advances in group theory and low-dimensional topology
  (Pusan, 2000)}, volume~27 of {\em Res. Exp. Math.}, pages 1--9. Heldermann,
  Lemgo, 2003.

\bibitem{B}
A.~Borel.
\newblock Commensurability classes and volumes of hyperbolic {$3$}-manifolds.
\newblock {\em Ann. Scuola Norm. Sup. Pisa Cl. Sci. (4)}, 8(1):1--33, 1981.

\bibitem{CF}
T.~Chinburg and E.~Friedman.
\newblock An embedding theorem for quaternion algebras.
\newblock {\em J. London Math. Soc. (2)}, 60(1):33--44, 1999.

\bibitem{Co}
H.~Cohen~et. al.
\newblock List of number fields of degree $\leq 7$.
\newblock {\em available at ftp://math.u-bordeaux.fr/pub/number fields}.

\bibitem{Co1}
H.~Cohen~et. al.
\newblock Pari.
\newblock {\em Freeware available by anonymous FTP from megrez@math.u-bordeaux,
  directory pub/pari}.

\bibitem{EEK}
A.~L. Edmonds, J.~H. Ewing, and R.~S. Kulkarni.
\newblock Torsion free subgroups of {F}uchsian groups and tessellations of
  surfaces.
\newblock {\em Invent. Math.}, 69(3):331--346, 1982.

\bibitem{Gr}
L.~Greenberg.
\newblock Maximal {F}uchsian groups.
\newblock {\em Bull. Amer. Math. Soc.}, 69:569--573, 1963.

\bibitem{J1}
S.~Johansson.
\newblock Genera of arithmetic {F}uchsian groups.
\newblock {\em Acta Arith.}, 86(2):171--191, 1998.

\bibitem{J2}
S.~Johansson.
\newblock On fundamental domains of arithmetic {F}uchsian groups.
\newblock {\em Math. Comp.}, 69(229):339--349, 2000.

\bibitem{K}
S.~Katok.
\newblock {\em Fuchsian groups}.
\newblock Chicago Lectures in Mathematics. University of Chicago Press,
  Chicago, IL, 1992.

\bibitem{M}
M.~L. Macasieb.
\newblock Derived arithmetic {F}uchsian groups of genus two.
\newblock {\em Ph. D. Thesis}, 2005.

\bibitem{MRe}
C.~Maclachlan and A.~W. Reid.
\newblock {\em The arithmetic of hyperbolic 3-manifolds}, volume 219 of {\em
  Graduate Texts in Mathematics}.
\newblock Springer-Verlag, New York, 2003.

\bibitem{MR1}
C.~Maclachlan and G.~Rosenberger.
\newblock Two-generator arithmetic {F}uchsian groups.
\newblock {\em Math. Proc. Cambridge Philos. Soc.}, 93(3):383--391, 1983.

\bibitem{MR2}
C.~Maclachlan and G.~Rosenberger.
\newblock Two-generator arithmetic {F}uchsian groups. {II}.
\newblock {\em Math. Proc. Cambridge Philos. Soc.}, 111(1):7--24, 1992.

\bibitem{O}
A.~M. Odlyzko.
\newblock Some analytic estimates of class numbers and discriminants.
\newblock {\em Invent. Math.}, 29(3):275--286, 1975.

\bibitem{Pa}
C.~J. Parry.
\newblock Units of algebraic number fields.
\newblock {\em J. Number Theory}, 7(4):385--388, 1975.

\bibitem{octic}
M.~Pohst, J.~Martinet, and F.~Diaz~y Diaz.
\newblock The minimum discriminant of totally real octic fields.
\newblock {\em J. Number Theory}, 36(2):145--159, 1990.

\bibitem{Rib}
P.~Ribenboim.
\newblock {\em Algebraic numbers}.
\newblock Wiley-Interscience [A Division of John Wiley\thinspace \&\thinspace
  Sons, Inc.], New York-London-Sydney, 1972.
\newblock Pure and Applied Mathematics, Vol. 27.

\bibitem{Sc}
V.~Schneider.
\newblock Die elliptischen {F}ixpunkte zu {M}odulgruppen in
  {Q}uaternionenschiefk\"orpern.
\newblock {\em Math. Ann.}, 217(1):29--45, 1975.

\bibitem{T1}
K.~Takeuchi.
\newblock Commensurability classes of arithmetic triangle groups.
\newblock {\em J. Fac. Sci. Univ. Tokyo Sect. IA Math.}, 24(1):201--212, 1977.

\bibitem{T2}
K.~Takeuchi.
\newblock Arithmetic {F}uchsian groups with signature {$(1;e)$}.
\newblock {\em J. Math. Soc. Japan}, 35(3):381--407, 1983.

\bibitem{V}
M.-F. Vign{\'e}ras.
\newblock {\em Arithm\'etique des alg\`ebres de quaternions}, volume 800 of
  {\em Lecture Notes in Mathematics}.
\newblock Springer, Berlin, 1980.

\bibitem{W}
L.~C. Washington.
\newblock {\em Introduction to cyclotomic fields}, volume~83 of {\em Graduate
  Texts in Mathematics}.
\newblock Springer-Verlag, New York, second edition, 1997.

\end{thebibliography}

      % is inserted.                         %
\index{Bibliography@\emph{Bibliography}}%  

\begin{center}
\noindent\rule{4cm}{.5pt}
\vspace{.25cm}
\flushleft
\noindent {\sc \small Melissa L.~Macasieb}\\
{\small Department of Mathematics \\ University of British Columbia\\ Vancouver BC V6T 1Z2} \\
email: {\tt macasieb@math.ubc.ca}
\end{center}

\end{document}